%% file: arxiv_version.tex
\documentclass[12pt, a4paper, english]{amsart}
\usepackage{amsmath} 
\usepackage{amsthm}
\usepackage{amssymb}

\usepackage{csquotes}
\usepackage{svg}

\usepackage{lipsum} 
\usepackage{multicol}

\usepackage{amscd} 
\usepackage{mathtools}
\usepackage{xcolor}

\usepackage{array} 
\usepackage[backref=page,linktocpage]{hyperref} 
\usepackage{caption} 
\usepackage{graphics,graphicx} 
\usepackage{tikz,tikz-cd} 
\usetikzlibrary{graphs,graphs.standard,calc}
\usetikzlibrary{decorations.pathreplacing,}
\usetikzlibrary{automata}
\usetikzlibrary{fadings}
\usepackage{pgfplots}
\usepackage[utf8]{inputenc}

\usepackage{makecell}

\usepackage{booktabs}
\usepackage{xspace}
\newcommand{\TOPCOM}{\texttt{TOPCOM}\xspace}
\newcommand{\nauty}{\texttt{nauty}\xspace}
\newcommand{\mptopcom}{\texttt{mptopcom}\xspace}
\newcommand{\polymake}{\texttt{polymake}\xspace}

\usepackage{enumerate} 
\usepackage{newtxtext,newtxmath} 

\usepackage[margin=1.25in, vmargin=1.18in]{geometry}

\usepackage{verbatim}

\usepackage{scalerel}
\usepackage{standalone}

\DeclareMathOperator{\conv}{conv}

\hypersetup{
    colorlinks = true,
    linkbordercolor = {white},
    linkcolor = {purple},
    anchorcolor = {black},
    citecolor = {purple},
    filecolor = {cyan},
    menucolor = {red},
    runcolor = {cyan},
    urlcolor = {black}
}

\newtheoremstyle{theorems}
{13pt}
{13pt}
{\slshape}
{}
{\bfseries}
{}
{.5em}
{}

\theoremstyle{theorems}
\newtheorem{theorem}{Theorem}[section]
\newtheorem{corollary}[theorem]{Corollary}

\newtheorem{proposition}[theorem]{Proposition}
\newtheorem*{theorem-no-label}{Theorem}

\newtheoremstyle{definition}
{12pt}
{12pt}
{}
{}
{\bfseries}
{}
{.5em}
{}

\theoremstyle{definition}

\newtheorem{example}[theorem]{Example}
\newtheorem{remark}[theorem]{Remark}

\newtheorem{manualtheoreminner}{Theorem}
\newenvironment{manualtheorem}[1]{
  \IfBlankTF{#1}
    {\renewcommand{\themanualtheoreminner}{\unskip}}
    {\renewcommand\themanualtheoreminner{#1}}
  \manualtheoreminner
}{\endmanualtheoreminner}

\newcommand{\RR}{\mathbb{R}}
\newcommand{\CC}{\mathbb{C}}
\newcommand{\QQ}{\mathbb{Q}}
\newcommand{\ZZ}{\mathbb{Z}}
\newcommand{\PP}{\mathbb{P}}
\newcommand{\trop}{\text{trop}}
\newcommand{\val}{\text{val}}

\title{Tropical Elliptic Curves in 3-Space}
\subjclass[2020]{}

\author[L.~Casabella]{Laura Casabella}
\author[L.~Kastner]{Lars Kastner}
\author[R.~Vlad]{Raluca Vlad}

\address[Laura Casabella]{
	Max-Planck Institute for Mathematics in the Sciences, Leipzig, Germany}
\email{\url{laura.casabella@mis.mpg.de}}
\address[Lars Kastner]{
	Chair of Discrete Mathematics/Geometry, TU Berlin, Germany}
\email{\url{lkastner@math.tu-berlin.de}}
\address[Raluca Vlad]{
	Department of Mathematics, Brown University, Providence, RI 0291, USA}
\email{\url{raluca_vlad@brown.edu}}

\renewcommand*{\arraystretch}{1.1}
\newcommand*{\mline}[1]{
\begingroup
    \renewcommand*{\arraystretch}{1.1}
   \begin{tabular}[c]{@{}>{\centering\arraybackslash}p{2cm}@{}}#1\end{tabular}
  \endgroup
}
\usepackage{tabularray}

\begin{document}

\begin{abstract}
    We classify trivalent graphs with 16 vertices and 16 edges that	arise from
intersecting two quadratic surfaces in tropical 3-space. There are 4,009
such graphs, representing maximally degenerate stable models
of elliptic curves realized as tropical complete intersections of two quadrics. Our classification is derived from 405,246,030
regular unimodular triangulations of the 4-dimensional Cayley polytope.
\end{abstract}

\maketitle

\vspace{-5mm}

\section{Introduction}

Consider two quadratic polynomials $f_1,f_2 \in K[x,y,z]$ over some field $K$ with a non-archimedean valuation. 
The homogenizations of $f_1$ and $f_2$ define two quadratic surfaces in $\PP^3_K$, whose complete intersection is generically an elliptic curve.
We are interested in studying faithful tropicalizations of curves arising this way, as in Table~\ref{table:cycle-lengths}.

The Cayley polytope associated to our pair of polynomials is the $4$-dimensional tetrahedral prism $2\Delta_3 \times \Delta_1$, where the tetrahedron
$$2\Delta_3 := \conv\big\{\mathbf{0}, \, 2 \mathbf{e_1}, \, 2 \mathbf{e_2}, \, 2 \mathbf{e_3}\big\}  \; \subset \; \RR^3$$
is the Newton polytope of a quadratic surface, and $\Delta_1$ is the unit segment $[0,1] \subset \RR$. We will denote this Cayley polytope by $C(2\Delta_3, 2\Delta_3)$.

The coefficients of the two polynomials $f_1,f_2$ induce a regular polyhedral subdivision of $C(2\Delta_3, 2\Delta_3)$.
The combinatorial process of deriving the corresponding tropical elliptic curve from this subdivision relies on the theory of \emph{mixed subdivisions} and is summarized in Section~\ref{section:complete-int}, following~\cite[Section 4.6]{sturmfels-maclagan}.

In this paper, we are interested in \emph{smooth} tropical curves arising this way -- these are precisely those curves associated to \emph{unimodular triangulations} of the Cayley polytope. Each such smooth tropical curve is a trivalent combinatorial graph of genus $1$ with $16$ vertices, $16$ edges, and $16$ unbounded rays. We find the following.

\begin{theorem}\label{thm:main-thm-3-space}
    The Cayley polytope $C(2\Delta_3, 2\Delta_3)$ has 405,246,030 regular unimodular triangulations, up to the symmetry induced by the automorphism group $S_4 \times\ZZ_2$ of the polytope. These triangulations give rise to 4,009 distinct isomorphism classes of tropical elliptic curves.
\end{theorem}

The curves from our theorem form a complete list of tropical curves arising from two quadratic polynomials as described in~\cite[Proposition~4.6.20]{sturmfels-maclagan}. 
A consequence of our classification is that all possible cycle lengths are achieved (see Corollary~\ref{cor:all-cycle-lengths}).

The result was found
using the software \mptopcom~\cite{mptopcom-paper} and \texttt{polymake}~\cite{DMV:polymake}.
We provide a list of all of the triangulations and the isomorphism classes of graphs from Theorem~\ref{thm:main-thm-3-space}, together with the code we used, at the following link:
\begin{equation}\label{eq:link}
\text{\url{https://github.com/dmg-lab/tropical_elliptic_curves}.}
\end{equation}
We refer the reader to Sections~\ref{sec:tropical-elliptic-curves-3-space} and~\ref{section:computations} for computational details. See Remark~\ref{rem:entry-in-list-of-4009-graphs} for a summary of the data included in our atlas of 4,009 tropical elliptic curves in $3$-space.

\smallskip

On the one hand, the importance of our result lies in the significant computational challenge posed by the problem of generating all of the triangulations from Theorem~\ref{thm:main-thm-3-space} and, then, exploiting this data to obtain the graph isomorphism classes. On the other hand, our result also has geometric interpretations, which we now explain.

Every tropical elliptic curve $G$ that we obtain represents a faithful (i.e.\ revealing the full genus) tropicalization of an elliptic curve realized as the intersection of two quadrics in $\PP^3_K$.
The $16$ rays of $G$ correspond to the $16$ marked points on the classical elliptic curve obtained by intersecting the curve with the four coordinate planes in $\PP^3$. 
The tropical curve $G$ is therefore dual to a stable model of an elliptic curve with $16$ marked points. We refer to~\cite{Chan-lecture-notes} for details on this semistability perspective
of tropicalization.

Since the graph $G$ is trivalent, 
the corresponding stable curve is maximally degenerate. Equivalently, the graph $G$ encodes a deepest boundary stratum in the Deligne-Knudsen-Mumford compactification $\overline{M}_{1,16}$ of the moduli space of smooth curves with genus $1$ and $16$ marked points~\cite{deligne-mumford, knudsen}. These deepest boundary strata in the moduli space give rise to the most combinatorially rich data. Therefore, our 4,009 graphs from Theorem~\ref{thm:main-thm-3-space} provide a complete list of maximally degenerate stable models of elliptic curves dual to tropical complete intersections of two quadrics in projective $3$-space.

\medskip

\noindent \textbf{Further directions.} Our result suggests new directions of study in the theory of tropical curves. For example, one could ask for a similar classification of tropical complete intersections of surfaces of higher degrees in $\PP^3$. Understanding such tropical curves of higher genus would require different techniques, as our algorithms would not terminate given that the amount of data (e.g.\ the number of triangulations of the appropriate Cayley polytope) increases very fast with increasing degrees.
In this sense, our algorithms really pushed current computational capabilities close to the farthest extent possible in terms of software and hardware.

In this paper, we report on combinatorial types of graphs. However, tropical curves have the additional structure of \emph{metric} graphs -- i.e.\ their edges carry lengths, as described in~\cite{moduli-tropical-plane-curves, Chan-lecture-notes}.
It would be an interesting further direction to understand the locus of metric graphs that gets realized by tropicalizing complete intersections in $3$-space.

\medskip

\noindent \textbf{Related works.} 
In~\cite{moduli-tropical-plane-curves}, the authors study the space of metric graphs realized as tropical plane curves of genus at least $2$.
Tropical elliptic curves expressed as tropical plane cubics have been studied in~\cite{j-invariant} via the tropical $j$-invariant and in~\cite{honeycomb-form} via honeycomb triangulations.
In~\cite{CGP-marked-points}, the authors study the moduli space of tropical curves of genus $1$ with $n$ marked points, viewed as abstract metric graphs, without reference to an embedding. 
Faithful embeddings of hyperelliptic curves have been studied in~\cite{cueto-markwig} for genus $2$ and in~\cite{Markwig-Ristau-Schleis} for genus at least $3$.

On the computational side, the algorithm we use for listing regular triangulations 
was originally implemented in \cite{mptopcom-paper}, where the authors develop \mptopcom to enumerate 
tropical hypersurfaces.
The algorithm was subsequently applied to cyclic polytopes in \cite{cyclic}. 
It was further improved in \cite{mptopcom-flips}, and employed in \cite{Casabella12112024} to enumerate regular triangulations of hypersimplices. To our knowledge, the computation in our paper is the largest of its~kind. 

\medskip

\noindent \textbf{Organization.} The remainder of the paper proceeds as follows. In Section~\ref{sec:planar-case}, we classify tropical elliptic \emph{plane} curves (see Figure~\ref{fig:all-graphs}), as a warm-up to our main results in $3$-space. Section~\ref{section:complete-int} summarizes the process of tropicalizing generic complete intersections of two surfaces in $3$-space. Finally, Section~\ref{sec:tropical-elliptic-curves-3-space} presents our main result (Theorem~\ref{thm:main-thm-3-space}) and data, while Section~\ref{section:computations} is dedicated to computational details.

\medskip

\noindent \textbf{Acknowledgements.} The authors are grateful to Bernd Sturmfels for suggesting this problem, to Maximilian Wiesmann for the computation in Remark \ref{rem:total-graphs}, and to Shelby Cox, Michael Joswig, Marta Panizzut, and Smita Rajan for helpful comments. The second author was supported by \href{https://www.mardi4nfdi.de/}{MaRDI (Mathematical Research
Data Initiative)}, funded by the Deutsche Forschungsgemeinschaft (DFG), project
number 460135501, NFDI 29/1 ``MaRDI -- Mathematische
Forschungsdateninitiative.''
The third author thanks the Max Planck Institute for Mathematics in the Sciences for the 
hospitality
that made this project possible.

\section{The planar case}
\label{sec:planar-case}

We first describe, in Proposition~\ref{prop:planar-case} and Figure~\ref{fig:all-graphs}, the graphs underlying smooth cubic curves in the tropical plane. Our enumeration provides a complete combinatorial classification of the tropical elliptic curves from Exercise 13 in~\cite[Section~1.9]{sturmfels-maclagan}. 

Even though this $2$-dimensional case is a straightforward exercise in planar tropical geometry, we include it here because we are not aware of this classification explicitly appearing in the literature, and because it is a natural point of departure for our paper. 
Later, we will generalize this plane cubic case to the setting of elliptic curves embedded in higher dimensional tropical space. 
The planar results of this section also highlight, by comparison, the computational difficulty of our main results in $3$-space (see Remark~\ref{rem:computational-challenge-3-space}).

\begin{remark}
    Another natural generalization of the planar case from this section is that of higher dimensional tropical hypersurfaces. In tropical $3$-space, cubic surfaces have been studied in~\cite{schlafi}, while quartic surfaces and their K3 polytopes have been studied in~\cite{k3-polytopes}. This direction of study will not be of further concern in this paper.
\end{remark}

Consider thus a plane cubic polynomial
$$f \; = \sum_{(i,j)\in P \cap \ZZ^2} c_{i,j} \cdot x^iy^j \;\; \in \; K[x,y], \quad \text{ with all }\; c_{i,j} \neq 0,$$
where $P = 3\Delta_2 := \mathrm{conv}\big\{\mathbf{0}, \, 3\mathbf{e_1}, \, 3\mathbf{e_2}\big\}  \subset  \RR^2$ is its Newton polygon.
Let $T$ be the regular subdivision of $P$ induced by the valuations of the coefficients $c_{i,j}$, as explained in~\cite[Section~1.3]{sturmfels-maclagan}.
And let $G$ be the dual graph to the subdivision $T$, having one vertex for each maximal (i.e.\ $2$-dimensional) cell of $T$, one edge for each pair of maximal cells sharing a common facet, and one unbounded ray for each exterior edge of $T$.

We make the technical assumption that $T$ is a unimodular triangulation, which means that
$G$ is a \emph{smooth} tropical curve, in the sense of~\cite[Section~1.3]{sturmfels-maclagan}. The polygon $P$ has one interior lattice point, so the tropical curve $G$ is indeed elliptic.

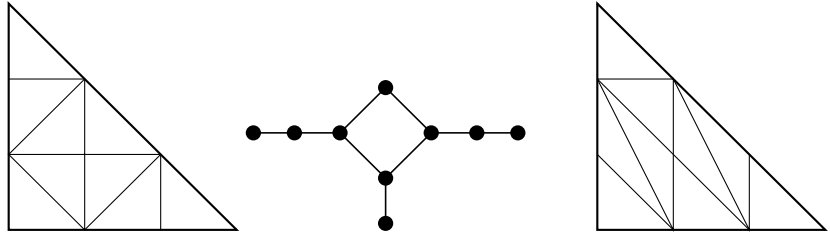
\begin{figure}[h]
\begin{center}
\begin{tikzpicture}
  \draw[thick] (0,0) -- (3,0) -- (0,3) -- cycle;

  \draw (0,2) -- (1,2) -- (1,1) -- (0,1) -- (1,0) -- (1,1) -- (2,1) -- (2,0);
  \draw (0,1) -- (1,2);
  \draw (1,0) -- (2,1);
\end{tikzpicture}
\scalebox{1.5}{\input{434957730.tikz}} \hspace{16pt}
\begin{tikzpicture}
  \draw[thick] (0,0) -- (3,0) -- (0,3) -- cycle;
  \draw ((1,2) -- (0,2) -- (1,0) -- (1,2) -- (2,0) -- (2,1);
  \draw (0,2) -- (2,0);
  \draw (0,1) -- (1,0);
\end{tikzpicture}
\end{center}
\caption{A tropical elliptic curve (depicted with rays suppressed) dual to two 
regular unimodular triangulations of $P = 3\Delta_2$.}
\label{fig:graph}
\end{figure}

\begin{figure}[h]
\begin{center}
\section*{Cycle length: 3}	
\input{1435761386.tikz} \qquad
\input{1163550515.tikz}
\section*{Cycle length: 4}
\input{1702484648.tikz}
\input{24790958.tikz}
\input{434957730.tikz}
\input{293328903.tikz}
\section*{Cycle length: 5}
\input{27921636.tikz} \quad
\input{552447728.tikz} \quad
\input{1837752289.tikz}	\quad
\input{562584937.tikz} 
\section*{Cycle length: 6}
\input{827130457.tikz} \quad
\input{877228084.tikz} \quad
\input{1095260065.tikz} \quad
\input{1347018042.tikz}

\begin{multicols}{3}
\section*{Cycle length: 7}
\input{2128804136.tikz}  \quad
\input{967574001.tikz}
\section*{Cycle length: 8}
\vspace{0.2cm}
\input{1623687162.tikz}
\section*{Cycle length: 9}
\vspace{0.2cm}
\input{1341331914.tikz}
\end{multicols}

\end{center}
    \caption{Isomorphism classes of smooth tropical elliptic plane curves, ordered in increasing order of cycle length.}
    \label{fig:all-graphs}
\end{figure}
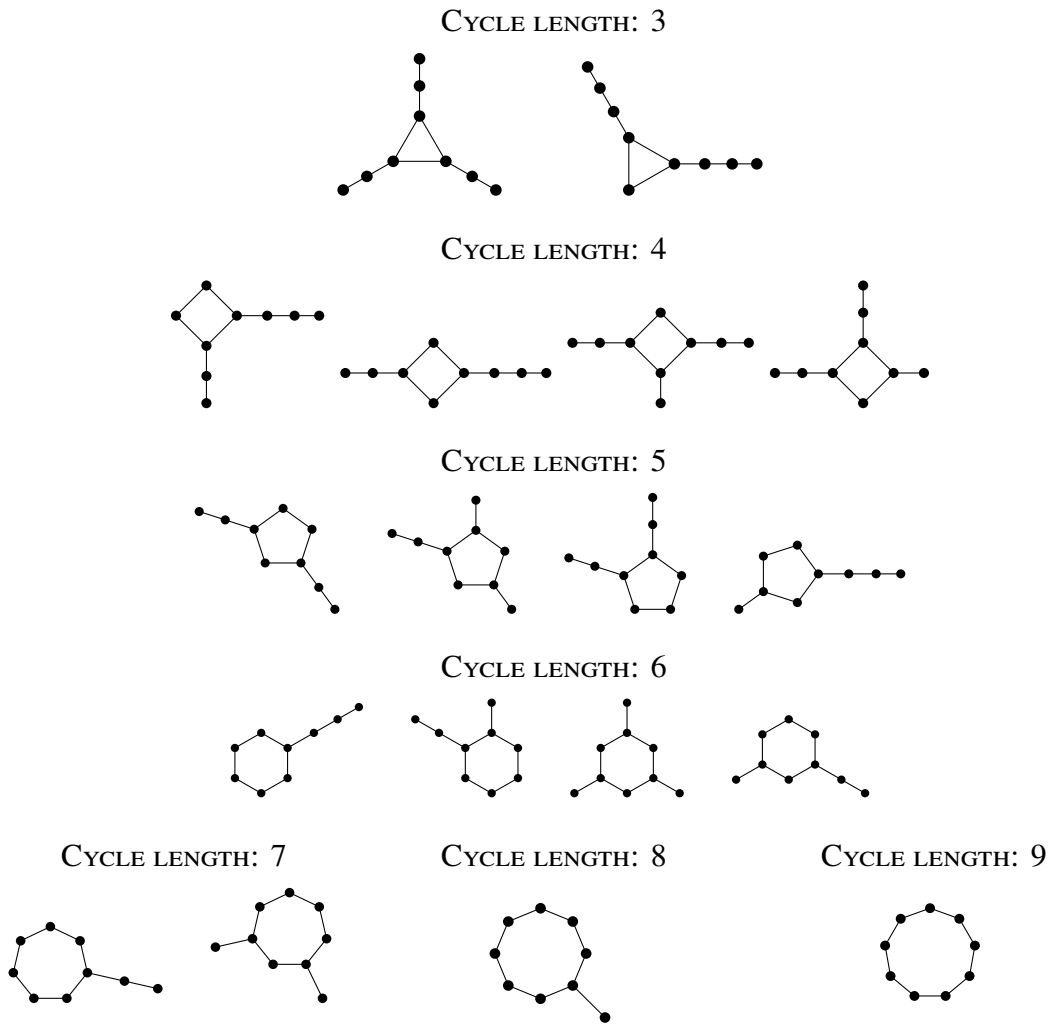

The tropical elliptic curve $G$ is a trivalent graph of genus $1$ with $9$ vertices, $9$ edges, and $9$ unbounded rays. We will suppress unbounded rays in this paper on account that their placement along $G$ is uniquely determined by the condition that the resulting graph should be trivalent. See Figure~\ref{fig:graph}.

Using \texttt{Gfan} \cite{gfan} through its interface in \texttt{polymake}, we generated all regular triangulations of $P$ and then filtered them by unimodularity. This complete enumeration relies on a computation of the secondary fan of $P$, which is a complete fan in $\RR^{10} = \RR^{\#(P \cap \ZZ^2)}$ encoding all regular subdivisions of $P$ via the weight vectors inducing them, see \cite[Section 5.2]{triangulations}. Using \texttt{polymake}, we then generated the dual graph for each unimodular triangulation and sorted these graphs into isomorphism classes, obtaining the following.

\begin{proposition}\label{prop:planar-case}
    The polygon $P = 3\Delta_2$ has 79 distinct regular unimodular triangulations,
    giving rise to 18 distinct isomorphism classes of tropical elliptic plane curves.
\end{proposition}

The 18 isomorphism classes of graphs are shown in Figure~\ref{fig:all-graphs}.

\begin{remark}
    The automorphism group $S_3$ acts on the triangle $P = 3\Delta_2$ by permuting its vertices. Under this symmetry action, the triangulations from Proposition~\ref{prop:planar-case} form $18$ orbits of regular unimodular triangulations of $P$. These symmetry classes of triangulations are therefore in bijection with the isomorphism classes of their dual graphs. For example, the two triangulations in Figure~\ref{fig:graph} are the same up to a rotation of $P$. This phenomenon will not persist in $3$-space, where different symmetry classes of triangulations can give rise to isomorphic tropical curves (see Theorem~\ref{thm:main-thm-3-space}).
\end{remark}

\begin{remark}
    The total number of connected trivalent graphs of genus $1$, having $9$ vertices and $9$ unbounded rays is $80$. Therefore, $22.5\%$ of these graphs are realized as tropical elliptic plane curves as in Proposition~\ref{prop:planar-case}. See Remark~\ref{rem:total-graphs} for details on the computation of the total number of graphs and an analogous percentage in $3$-space.
\end{remark}

\begin{remark}
    The tropical curve $G$ is a faithful tropicalization of the elliptic curve $V(\Tilde{f}) \subset \PP_K^2$, where $\Tilde{f}$ denotes the homogenization of the polynomial $f$. 
    Analogous to the geometric interpretation put forth in the Introduction, the tropical plane curve $G$ encodes a maximally degenerate stable model of a $9$-marked elliptic curve. On the algebraic side, the $9$ marked points in question are the intersection of the plane cubic curve $V(\Tilde{f})$ with the three coordinate lines of $\PP_K^2$. On the tropical side, the marked points correspond to the unbounded rays.
    At the moduli level, 
    our tropical curve $G$ encodes a deepest boundary stratum in the Deligne-Knudsen-Mumford compactification $\overline{M}_{1,9}$ of the moduli space of smooth genus $1$ curves with $9$ marked points~\cite{deligne-mumford, knudsen}. 
\end{remark}

\section{Tropicalizations of complete intersections} \label{section:complete-int}

In this section, we summarize the construction of the tropicalization of a complete intersection of two surfaces in $3$-space. We refer the reader to~\cite[Section~4.6]{sturmfels-maclagan} for further details and analogous procedures for tropicalizations of more general complete intersections. In Section~\ref{sec:tropical-elliptic-curves-3-space}, we will specialize this construction to elliptic curves realized as complete intersections of two quadratic surfaces. 

\smallskip

Let $\Delta_3 := \mathrm{conv}\big\{ \mathbf{0}, \mathbf{e_1}, \mathbf{e_2}, \mathbf{e_3}\big\} \subset \RR^3$ denote the standard tetrahedron. Let $f_1,f_2 \in K[x,y,z]$ be two polynomials of degree $d$ and $e$, respectively, having full monomial support. Their Newton polytopes are $P_1 = d \Delta_3$ and  $P_2 = e\Delta_3$, respectively, where $n \Delta_3$ denotes the dilation of the standard simplex $\Delta_3$ by a factor of $n$.

The \emph{Cayley polytope} of $P_1$ and $P_2$ is the tetrahedral prism given by
\begin{equation}\label{eq:cayley-polytope}
    C(P_1, P_2) \coloneqq \conv\big\{ \mathbf{e_1} \times P_1, \mathbf{e_2} \times P_2\big\} \; \subset \;  \RR^{2} \times \RR^3 = \RR^5. 
\end{equation}

The \emph{Minkowski sum} $P_1+P_2 \coloneqq \{p_1 + p_2 \mid p_i \in P_i\} \subset \RR^3$ of the two Newton polytopes is isomorphic, via scaling by a factor of $1/2$, to the following slice of the Cayley polytope:
\begin{equation}\label{eq:minkowski-cayley-slice}
    P_1 + P_2 \; \cong \; C(P_1,P_2) \cap \{x_1 = x_2 = 1/2\}, 
\end{equation}
where $x_1,x_2$ denote the coordinates on the $\RR^2$ factor of the ambient $\RR^5$ in~\eqref{eq:cayley-polytope}.

\smallskip

There is a one-to-one correspondence between the set of lattice points of $C(P_1,P_2)$ and the disjoint union of the sets of lattice points of $P_1$ and $P_2$. We are assuming that the valuations of the coefficients of the polynomials $f_1, f_2$ are generic, so they induce a regular polyhedral triangulation of 
the Cayley polytope $C(P_1,P_2)$.
We are again primarily interested in the case where the induced triangulation 
is unimodular. 

\begin{remark}
    By Sturmfels' correspondence~\cite[Theorem 9.4.5]{triangulations}, our unimodularity condition is equivalent to the initial ideal of the toric ideal $I_{C(P_1,P_2)}$ being square-free, where the initial ideal is taken with respect to the term order induced by the valuations of the coefficients of the polynomials $f_1, f_2$.
\end{remark}

The associated \emph{smooth} tropical curve
$$G \; := \;\trop\big(V(f_1,f_2)\big) \, = \; \trop\big(V(f_1)\big) \cap \, \trop\big(V(f_2)\big)$$
is a faithful tropicalization of the complete intersection $V(f_1,f_2) \subset (K^*)^3$. Moreover, by~\cite[Theorem 4.6.18]{sturmfels-maclagan}, this tropical curve can be obtained from the following combinatorial process. Our triangulation of $C(P_1,P_2)$ induces a \emph{mixed subdivision} of the Minkowski sum $P_1 + P_2$ via (\ref{eq:minkowski-cayley-slice}). The maximal cells in the latter subdivision are of the form
$Q_1 + Q_2,$ 
where $Q_i \subset P_i$ is a cell in the polyhedral subdivision of $P_i$ induced by the valuations of the coefficients of $f_i$, and $\dim(Q_1) + \dim(Q_2) = 3$. 
We are interested in the \emph{mixed cells} where $\dim(Q_i) > 0$ for both $i = 1,2$. In our setting, the mixed cells are always \enquote{\emph{toblerones}} (i.e.\ triangular prisms); such a toblerone is obtained as the Minkowski sum $Q_1 + Q_2$ where one cell $Q_i$ is a segment and the other is a triangle.

Then, our tropical curve $G$ consists of one vertex for every mixed cell, with an edge between any two vertices corresponding to adjacent toblerones sharing a quadrilateral facet. We will again suppress the rays usually associated to mixed cells having exterior quadrilateral facets; these rays attach to the graph $G$ in the unique way that makes the resulting graph trivalent.

\medskip

To end this section, we illustrate this mixed subdivision construction for an example in degrees $d = 2$ and $e = 1$. Even though our main results concern complete intersections of two surfaces of degree $2$ (see Sections~\ref{sec:tropical-elliptic-curves-3-space} and~\ref{section:computations}), we depict a smaller example here in order to make the figures intelligible, with all mixed cells visible.

\begin{figure}[h]
    \centering
\input{minkowski3.tex}
    \caption{The mixed subdivision of $3\Delta_3$ from Example~\ref{ex:mixed-subdivision} with the toblerones highlighted in blue (left), the toblerones of this subdivision isolated for clarity (middle), and the dual smooth tropical curve (right).}
    \label{fig:toblerones}
\end{figure}
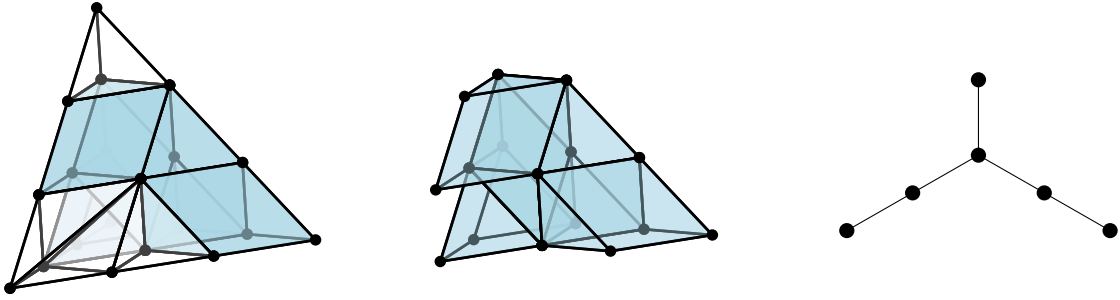

\begin{example}\label{ex:mixed-subdivision}
Consider the field of complex Puiseux series $K = \CC\{\!\{t\}\!\}$ with the valuation $\val: K \to \QQ \cup \{\infty\}$ sending a non-zero element $\sum_{m = m_0}^\infty c_mt^{m/n}$, with $c_{m_0} \neq 0$, to $m_0/n$, which is the smallest power of $t$ appearing in the series; by convention, $\val(0) = \infty$.

Let $f_1, f_2 \in K[x,y,z]$ be the polynomials given by 
\begin{align*}
    f_1 &= t^3x+xy+t^2xz+x+t^4y^2+t^2yz+y+t^3z^2+z+t^3\\
    \text{and } \;\;f_2 &= t^2x+ty+z+t^2,
\end{align*}
having Newton polytopes $P_1 = 2\Delta_3$ and $P_2 = \Delta_3$, respectively. The coefficients of the polynomials $f_1, f_2$ induce a unimodular triangulation of the associated Cayley polytope $C(2\Delta_3, \Delta_3)$. Slicing according to~\eqref{eq:minkowski-cayley-slice}, we obtain a subdivision of the Minkowski sum $P_1 + P_2 = 3\Delta_3$ into nine simplices and six toblerones. See Figure~\ref{fig:toblerones}.

\end{example}

\section{Tropical elliptic curves in 3-space} \label{sec:tropical-elliptic-curves-3-space}

We now present our main computation.
Consider thus generic polynomials $f_1, f_2 \in K[x,y,z]$ of degree~$2$. 
The tropicalization of the curve $Y = V(f_1,f_2) \subset (K^*)^3$ from Section~\ref{section:complete-int} is equivalently the tropicalization of the closure $\overline{Y} \subset \PP^3_K$ with the $16$ points $\overline{Y} \setminus Y$ being marked; here, $\overline{Y}$ is the elliptic curve cut out by the homogenizations of $f_1$ and $f_2$.
Elliptic curves realized as such complete intersections 
 of two quadrics in projective $3$-space 
have recent applications to game theory, as studied in~\cite{elke-irem}.

The Newton polytopes of $f_1,f_2$ are $P_1=P_2=2\Delta_3$. The associated Cayley polytope $C(P_1,P_2) = C(2\Delta_3, 2\Delta_3)$ is a four-dimensional tetrahedral prism in $\RR^5$ having $8$ vertices and $20$ lattice points. 
Each unimodular triangulation of this Cayley polytope induces a mixed subdivision on the Minkowski sum $P_1+P_2 = 4\Delta_3$, as described in Section~\ref{section:complete-int}. 

The Cayley polytope $C(2\Delta_3, 2\Delta_3)$ has lattice volume $32$, so a unimodular triangulation of this polytope consists of $32$ four-dimensional simplices. Therefore, the induced mixed subdivision of the Minkowski sum consists of $32$ maximal cells. Since the Minkowski sum $4\Delta_3$ has lattice volume $64$, and a triangular prism (or toblerone) has lattice volume $3$, a quick computation shows that 
the induced Minkowski subdivision must consist of precisely $16$ mixed cells (toblerones) and $16$ \emph{unmixed} cells (unit simplices obtained as $Q_1 + Q_2$ with one of the cells $Q_i$ being a point and the other a tetrahedron).

A tropical curve $G$ dual to such a mixed subdivision is a trivalent graph of genus $1$ with $16$ vertices (corresponding to the $16$ mixed cells), $16$ edges, and $16$ unbounded rays (which we suppress).
Our aim is to classify all combinatorial types of smooth tropical elliptic curves obtained from this combinatorial process. We find the following result, which was established through computational methods, with algorithms, software, and runtime documented in Section~\ref{section:computations}.

\begin{manualtheorem}{\ref{thm:main-thm-3-space}}
\emph{The Cayley polytope $C(2\Delta_3, 2\Delta_3)$ has $405,\!246,\!030$ regular unimodular triangulations, up to the symmetry induced by the automorphism group $S_4 \times\ZZ_2$ of the polytope. These triangulations give rise to $4,\!009$ distinct isomorphism classes of tropical elliptic curves.}
\end{manualtheorem}

The automorphism group of our Cayley polytope is the direct product $S_4\times\ZZ_2$. Here, the symmetric group $S_4$ acts simultaneously on the two copies of $2\Delta_3$ inside $C(2\Delta_3, 2\Delta_3)$, by permuting the vertices of each. Additionally, $\ZZ_2$ exchanges the first two coordinates in~\eqref{eq:cayley-polytope}.
This action is affine linear, so it induces an action on the set of all triangulations.

Our code, together with the data of all the triangulations and the isomorphism classes of graphs from Theorem~\ref{thm:main-thm-3-space}, can be downloaded at the link~\eqref{eq:link}. 
We also counted all (not necessarily regular) unimodular triangulations of $C(2\Delta_3, 2\Delta_3)$ -- there are 20,360,876,352 such triangulations, forming 424,262,728 symmetry classes under the $(S_4 \times \ZZ_2)$-action.
Our data also includes
an enumeration of the regular full (i.e.\ using all lattice points) triangulations of $C(2\Delta_3, 2\Delta_3)$ -- there are 836,381,618 such triangulations, up to symmetry. 
See Section~\ref{section:computations}
for further details and documentation.

\begin{remark}\label{rem:entry-in-list-of-4009-graphs}
    \textbf{(Database)} The $336$-page file \texttt{atlas.pdf} at the link~\eqref{eq:link} contains the atlas of the 4,009 isomorphism classes of graphs from Theorem~\ref{thm:main-thm-3-space}, listed in increasing order of cycle length. In Figure~\ref{fig:sample-entry}, we include one sample entry from the atlas and we explain here the information provided therein. 

    \begin{figure}[h]
    \centering
    \includegraphics[width=0.8\linewidth]{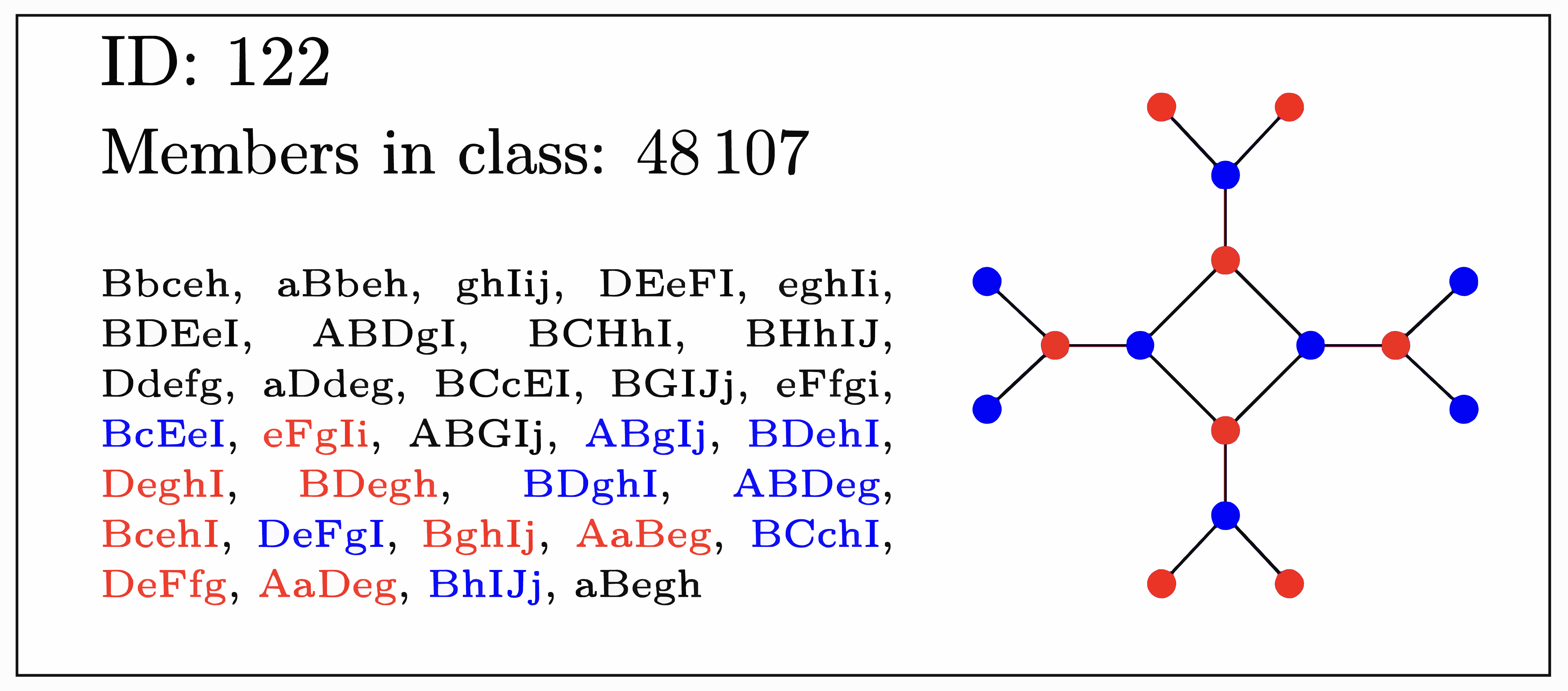}
    \caption{One of the 4,009 entries in our atlas of tropical elliptic curves.}
    \label{fig:sample-entry}
\end{figure}

    {\setlength{\leftmargini}{18pt}
    \begin{enumerate}
        \item The ID of the graph runs from $0$ to 4,008 and records its place in our list.

        \smallskip
        
        \item The number of ``members in class" records how many symmetry classes of unimodular triangulations of the Cayley polytope $C(2\Delta_3, 2\Delta_3)$ give rise to a tropical elliptic curve in that given graph isomorphism class.

        \smallskip
        
        \item For every graph, we also record one unimodular triangulation of $C(2\Delta_3, 2\Delta_3)$ inducing it. We explain here how this triangulation is encoded. 
        
        The Cayley polytope $C(2\Delta_3, 2\Delta_3) = \conv\big\{ \mathbf{e_1} \times 2\Delta_3, \mathbf{e_2} \times 2\Delta_3\big\}$ contains $20$ lattice points, half of them coming from $\mathbf{e_1} \times 2\Delta_3$ and the other half from $\mathbf{e_2} \times 2\Delta_3$. We label the first $10$ lattice points by $A,\dots,J$ and the other $10$ by $a,\dots,j$. The precise bijection between our labels and the lattice points of $C(2\Delta_3, 2\Delta_3)$ can be found on page $1$ of the file \texttt{atlas.pdf}.

        Any unimodular triangulation of $C(2\Delta_3, 2\Delta_3)$ consists of $32$ maximal cells, each of which is a four-dimensional simplex. The $5$ vertices of each simplex consist of $k$ lowercase and $5-k$ uppercase lattice points, for some $k \in \{1,2,3,4\}$. We write our sample triangulation via its 32 maximal cells, encoded through their corresponding 5-letter strings. For example, \enquote{BcEeI} is a four-dimensional simplex with $k = 2$.

        The toblerones in the induced Minkowski subdivision are precisely the slices of those simplices with $k = 2$ or $3$. We color the strings of these simplices in our sample triangulation -- the colors red and blue distinguish between the cases $k = 2$ and $3$. For example, the strings colored blue correspond to simplices in the Cayley polytope triangulation that have $2$ lowercase and $3$ uppercase vertices. These are precisely the simplices inducing toblerones $Q_1 + Q_2$ where $Q_1 \subset P_1 = 2\Delta_3$ is the triangle whose vertices are the $3$ uppercase points, and $Q_2 \subset P_2 = 2\Delta_3$ is the line segment whose endpoints are the $2$ lowercase points.

        The vertices of the tropical elliptic curve, which correspond to toblerones, are accordingly colored red/blue. We always have $8$ blue vertices and $8$ red vertices.
    \end{enumerate}}
\end{remark}

\begin{example}(\cite[Example~4.6.21]{sturmfels-maclagan})
    Consider $K = \QQ$ with the $2$-adic valuation, and
    \begin{align*}
        f_1 & = 1024 x^2 + 64xy + 8xz + 2x + 8y^2 + 2yz + y + z^2 + z + 2,\\
        f_2 &= x^2 + xy + 2xz + 8x + 2y^2 + 8yz + 64y + 64 z^2 + 1024z + 32768.
    \end{align*}
    The corresponding tropical elliptic curve has cycle length $8$ and has ID 2257 in our atlas.
\end{example}

\begin{remark}\label{rem:total-graphs}
    The total number of connected graphs of maximum degree 3 having genus 1, 16 vertices, and hence 16 edges, is 25,788. This number was obtained by computing the ordinary generating function of the combinatorial class of such graphs, using techniques described in \cite[Chapter I]{Flajolet_2009}. 
    Therefore, the percentage of such graphs realizable as smooth tropical elliptic curves from $C(2\Delta_3,2\Delta_3)$ is approx.\! $15.55\%$.
\end{remark}

\begin{remark}\label{rem:computational-challenge-3-space}
    Our computations highlight the difference in complexity
    between classifying tropical curves arising 
    from plane curves vs.\ from
    complete intersections of two surfaces in 3-space. The planar case of Section~\ref{sec:planar-case} could be handled in a few minutes using built-in functions from \texttt{polymake}, while the higher-dimensional case required several months of computations, the integration of multiple pieces of software, and access to a cluster. See Section~\ref{section:computations} for computational details. 
\end{remark}

\medskip

The following corollary is a consequence of our classification of tropical elliptic curves in $3$-space and answers the question which, in fact, originally motivated this~project.

\begin{corollary}\label{cor:all-cycle-lengths}
    Tropical elliptic curves from regular unimodular triangulations of the Cayley polytope $C(2\Delta_3, 2\Delta_3)$ achieve all possible cycle lengths (between $3$ and $16$).
\end{corollary}

In Table~\ref{table:cycle-lengths}, we provide a list of quadratic polynomials $f_1,f_2 \in \CC\{\!\{ t\}\!\}[x,y,z]$ inducing tropical elliptic curves of all possible cycle lengths.
For each cycle length, the 
coefficients of the corresponding pair of polynomials induce the sample triangulation associated to the first graph having that cycle length in our atlas.
We emphasize that one could work over any field with a non-trivial valuation; the tropical elliptic curves in Table~\ref{table:cycle-lengths} remain the same as long as the valuations of the coefficients are unchanged.

\begin{table}
\footnotesize
\begin{tblr}{colspec=|c|l|c|,row{1-14}={9ex}}
\hline
3 & \begin{tblr}{l} $t^{5}x^2 + xy + xz + x + t^{11}y^2 + t^{8}yz + t^{5}y + t^{9}z^2 + t^{3}z + t$ \\ $t^{4}x^2 + t^{2}xy + txz + x + t^{14}y^2 + yz + t^{9}y + t^{13}z^2 + t^{9}z + t^{6}$ \end{tblr} & \mline{\includestandalone{3}} \\ 
\hline 
4 & \begin{tblr}{l} $t^{22}x^2 + txy + t^{16}xz + t^{10}x + t^{8}y^2 + yz + t^{3}y + t^{11}z^2 + t^{2}z + 1$ \\ $t^{11}x^2 + xy + t^{4}xz + x + t^{19}y^2 + t^{7}yz + t^{15}y + z^2 + t^{6}z + t^{13}$ \end{tblr} & \mline{\includestandalone{4}} \\ 
\hline 
5 & \begin{tblr}{l} $t^{13}x^2 + t^{2}xy + t^{13}xz + t^{6}x + t^{13}y^2 + yz + t^{4}y + t^{14}z^2 + t^{3}z + 1$ \\ $tx^2 + t^{4}xy + xz + x + t^{16}y^2 + yz + t^{8}y + t^{6}z^2 + t^{2}z + t$ \end{tblr} & \mline{\includestandalone{5}} \\ 
\hline 
6 & \begin{tblr}{l} $t^{7}x^2 + xy + t^{5}xz + tx + t^{11}y^2 + t^{3}yz + t^{4}y + t^{7}z^2 + z + 1$ \\ $t^{5}x^2 + t^{14}xy + t^{2}xz + t^{15}x + t^{37}y^2 + yz + t^{31}y + t^{3}z^2 + z + t^{26}$ \end{tblr} & \mline{\includestandalone{6}} \\ 
\hline 
7 & \begin{tblr}{l} $t^{21}x^2 + t^{13}xy + t^{3}xz + t^{12}x + t^{12}y^2 + tyz + t^{8}y + z^2 + z + t^{6}$ \\ $t^{17}x^2 + xy + xz + t^{10}x + t^{3}y^2 + tyz + y + t^{2}z^2 + t^{3}z + t^{5}$ \end{tblr} & \mline{\includestandalone{7}} \\ 
\hline 
8 & \begin{tblr}{l} $t^{25}x^2 + t^{25}xy + t^{12}xz + t^{14}x + t^{26}y^2 + tyz + t^{12}y + z^2 + z + t^{4}$ \\ $tx^2 + xy + xz + t^{7}x + t^{50}y^2 + t^{24}yz + t^{32}y + z^2 + t^{5}z + t^{15}$ \end{tblr} & \mline{\includestandalone{8}} \\ 
\hline 
9 & \begin{tblr}{l} $t^{18}x^2 + t^{3}xy + t^{8}xz + t^{10}x + t^{4}y^2 + yz + t^{2}y + t^{6}z^2 + t^{2}z + t^{3}$ \\ $t^{9}x^2 + t^{2}xy + t^{2}xz + x + t^{4}y^2 + tyz + y + t^{5}z^2 + z + 1$ \end{tblr} & \mline{\includestandalone{9}} \\ 
\hline 
10 & \begin{tblr}{l} $t^{9}x^2 + t^{6}xy + xz + t^{4}x + t^{5}y^2 + t^{6}yz + y + t^{13}z^2 + z + 1$ \\ $t^{2}x^2 + xy + t^{5}xz + tx + t^{17}y^2 + t^{15}yz + t^{15}y + t^{17}z^2 + t^{14}z + t^{17}$ \end{tblr} & \mline{\includestandalone{10}} \\ 
\hline 
11 & \begin{tblr}{l} $t^{19}x^2 + t^{5}xy + t^{16}xz + t^{9}x + t^{2}y^2 + tyz + y + t^{15}z^2 + z + 1$ \\ $t^{13}x^2 + xy + t^{11}xz + t^{2}x + y^2 + tyz + t^{3}y + t^{11}z^2 + t^{8}z + t^{7}$ \end{tblr} & \mline{\includestandalone{11}} \\ 
\hline 
12 & \begin{tblr}{l} $t^{7}x^2 + t^{6}xy + t^{5}xz + t^{3}x + t^{15}y^2 + t^{12}yz + y + t^{10}z^2 + z + 1$ \\ $t^{9}x^2 + xy + t^{5}xz + t^{7}x + t^{4}y^2 + t^{6}yz + y + t^{9}z^2 + t^{7}z + t^{6}$ \end{tblr} & \mline{\includestandalone{12}} \\ 
\hline 
13 & \begin{tblr}{l} $t^{49}x^2 + t^{15}xy + t^{22}xz + t^{25}x + t^{2}y^2 + t^{3}yz + y + t^{9}z^2 + z + t^{2}$ \\ $t^{16}x^2 + xy + t^{5}xz + t^{8}x + ty^2 + yz + y + t^{4}z^2 + t^{2}z + t$ \end{tblr} & \mline{\includestandalone{13}} \\ 
\hline 
14 & \begin{tblr}{l} $t^{38}x^2 + t^{8}xy + t^{14}xz + t^{22}x + y^2 + t^{2}yz + y + t^{6}z^2 + z + t^{7}$ \\ $t^{14}x^2 + xy + t^{5}xz + t^{13}x + y^2 + tyz + t^{6}y + t^{4}z^2 + t^{5}z + t^{20}$ \end{tblr} & \mline{\includestandalone{14}} \\ 
\hline 
15 & \begin{tblr}{l} $t^{49}x^2 + t^{19}xy + t^{21}xz + t^{26}x + t^{8}y^2 + t^{5}yz + t^{4}y + t^{6}z^2 + z + t^{4}$ \\ $t^{14}x^2 + xy + t^{2}xz + t^{7}x + t^{5}y^2 + yz + t^{2}y + z^2 + z + t$ \end{tblr} & \mline{\includestandalone{15}} \\ 
\hline 
16 & \begin{tblr}{l} $t^{18}x^2 + xy + t^{7}xz + t^{9}x + t^{2}y^2 + t^{4}yz + y + t^{10}z^2 + z + t$ \\ $x^2 + xy + t^{7}xz + t^{9}x + t^{20}y^2 + t^{21}yz + t^{19}y + t^{26}z^2 + t^{22}z + t^{19}$ \end{tblr} & \mline{\includestandalone{16}} \\ 
\hline 
\end{tblr}
\caption{Pairs of quadratic polynomials in $\CC\{\!\{ t\}\!\} [x,y,z]$ inducing tropical elliptic curves of all possible cycle lengths (between $3$ and $16$).
} 
\label{table:cycle-lengths}
\end{table}

\begin{remark}
    Every elliptic curve can be written as a complete intersection of two quadrics in projective $3$-space, see \cite[Section 19.9.11]{vakil}. As before, we mark the $16$ points where the curve meets the coordinate planes of $\PP^3$. It remains an open question which are \emph{all} the tropical curves that are achieved as tropicalizations of such $16$-marked elliptic curves. For an unmarked example, see~\cite{tropicalizad-quartics} where the tropical curve of genus $3$ in Figure 13 cannot be achieved as a tropical plane quartic, but can be achieved as the tropicalization of a plane quartic after the re-embedding from Proposition~3.8.
\end{remark}

\section{Computational aspects} \label{section:computations}

In this section, we provide documentation for the algorithms we used. All our code and data can be accessed through the link~\eqref{eq:link}.
As before, when talking about triangulations of the Cayley polytope $C(2\Delta_3,2\Delta_3)$, we mean symmetry classes of triangulations under the $(S_4 \times \ZZ_2)$-action described in Section~\ref{sec:tropical-elliptic-curves-3-space}.

\medskip

\noindent \textbf{Computation strategy.} Our approach in computing the regular unimodular triangulations of $C(2\Delta_3,2\Delta_3)$ had two steps.

\setlength{\leftmargini}{18pt}{\begin{enumerate}
    \item We used \TOPCOM~\cite{TOPCOM} version 1.1.2 to count all unimodular triangulations of $C(2\Delta_3,2\Delta_3)$. This step, which was relatively fast because it did not check for regularity, was only meant as an initial test to see whether it would even be feasible to attempt a complete enumeration of all regular unimodular triangulations. 

We determined that $C(2\Delta_3,2\Delta_3)$ has 424,262,728 unimodular triangulations; the \TOPCOM computation took roughly 4 hours. Since this was only meant as an initial feasibility check, we did not collect the output, but instead used the \TOPCOM function
\texttt{points2nalltriangs} for only counting.

\smallskip

\item We generated all regular unimodular triangulations of $C(2\Delta_3,2\Delta_3)$ using
\mptopcom, which has a fast regularity check, see \cite{mptopcom-flips}.
We first generated all regular full triangulations of $C(2\Delta_3,2\Delta_3)$; there are 836,381,618 such triangulations. Afterwards, we ran a small script to filter out the non-unimodular ones and arrived at the 405,246,030 regular unimodular triangulations from Theorem~\ref{thm:main-thm-3-space}.
\end{enumerate}}

\medskip

\noindent \textbf{Resources.} We used the cluster of TU Berlin and allocated 30 slots to this computation. The underlying machines vary, but the operating system is AlmaLinux release 8.10 (Cerulean Leopard). An important feature of \mptopcom is checkpointing, i.e.\ it can halt the computation and write a checkpoint from which it may be restarted later. This is useful in several ways. On the one hand, large computing clusters do not usually allow jobs that run indefinitely long, so one has to specify a time limit. Furthermore, the scheduling algorithm on large clusters prioritizes jobs with shorter runtime as they are easier to squeeze into the schedule. 
Another aspect is that checkpointing prevents losing the entire computation in case of unexpected downtime. 
Having these advantages in mind, we chose the shortest possible time limit allowed by the cluster and ran our computation in chunks of 12 hours.
We ended up with 313 chunks of data. The total time for the computation was roughly 62 days.

\medskip

\noindent \textbf{Data processing.} We further processed our output using \polymake, as follows.

\setlength{\leftmargini}{18pt}{\begin{enumerate}
\item We extracted the regular unimodular triangulations. This was done using a small \texttt{perl} script and only took a few hours. 

\smallskip

\item For every triangulation, we computed the associated tropical curve, according to the construction described in Sections~\ref{section:complete-int} and~\ref{sec:tropical-elliptic-curves-3-space}, and determined its cycle length.

\smallskip

\item We identified the graph isomorphism classes of the generated tropical curves. For this, \polymake used \nauty \cite{nauty} in the background. To speed up the isomorphism class computation, we first computed what \nauty calls the \emph{canonical hash} of a graph. Essentially, two graphs cannot be isomorphic if their canonical hashes differ and, furthermore, the probability of hash collisions
(i.e.\ non-isomorphic graphs having the same hash)
is very low. We collected one representative for every canonical hash and, for every new graph, we only tested isomorphism if a representative had the same hash. We encountered no hash collisions.

\smallskip

\item Finally, we wrote code to turn the graphs, together with one sample triangulation per graph, into a nice visual atlas (see Remark~\ref{rem:entry-in-list-of-4009-graphs}).
\end{enumerate}}

\medskip

\noindent \textbf{Data.} The following data, uploaded on \texttt{GitHub} and \texttt{Zenodo}, can be accessed either through the link~\eqref{eq:link}, or directly under the following entries: 
\begin{itemize}
\item \cite{zenodo_regular_full} -- the 836,381,618 regular \emph{full} triangulations of $C(2\Delta_3,2\Delta_3)$;
\item \cite{zenodo_regular_unimodular} -- the 405,246,030 regular \emph{unimodular} triangulations of $C(2\Delta_3,2\Delta_3)$;
\item \texttt{atlas.pdf} in our repository~\eqref{eq:link} -- representatives for the graph isomorphism classes from Theorem~\ref{thm:main-thm-3-space}.
\end{itemize}

\bibliographystyle{alpha}
\bibliography{my}

\end{document}

%% file: 434957730.tikz
\begin{tikzpicture}[scale=0.4]
\draw (2.9,0) -- (2,0);
\draw (2,0) -- (1,0);
\draw (1,0) -- (6.12323399573677e-17,1);
\draw (1,0) -- (-1.83697019872103e-16,-1);
\draw (6.12323399573677e-17,1) -- (-1,1.22464679914735e-16);
\draw (-1.83697019872103e-16,-1) -- (-1,1.22464679914735e-16);
\draw (-1.83697019872103e-16,-1) -- (-3.67394039744206e-16,-2);
\draw (-1,1.22464679914735e-16) -- (-2,2.44929359829471e-16);
\draw (-2,2.44929359829471e-16) -- (-2.9,3.55147571752732e-16);
\fill (-1.83697019872103e-16, -1) circle (4.75pt);
\fill (2, 0) circle (4.75pt);
\fill (-2, 2.44929359829471e-16) circle (4.75pt);
\fill (-2.9, 3.55147571752732e-16) circle (4.75pt);
\fill (-1, 1.22464679914735e-16) circle (4.75pt);
\fill (2.9, 0) circle (4.75pt);
\fill (-3.67394039744206e-16, -2) circle (4.75pt);
\fill (1, 0) circle (4.75pt);
\fill (6.12323399573677e-17, 1) circle (4.75pt);
\end{tikzpicture}

%% file: 1435761386.tikz
\begin{tikzpicture}[scale=0.4]
\draw (0,2.9) -- (0,2);
\draw (0,2) -- (0,1);
\draw (0,1) -- (0.866025403784439,-0.5);
\draw (0,1) -- (-0.866025403784438,-0.5);
\draw (0.866025403784439,-0.5) -- (-0.866025403784438,-0.5);
\draw (0.866025403784439,-0.5) -- (1.73205080756888,-1);
\draw (-0.866025403784438,-0.5) -- (-1.73205080756888,-1);
\draw (1.73205080756888,-1) -- (2.51147367097487,-1.45);
\draw (-1.73205080756888,-1) -- (-2.51147367097487,-1.45);
\fill (2.51147367097487, -1.45) circle (5.43130215351701pt);
\fill (1.73205080756888, -1) circle (5.43130215351701pt);
\fill (0, 2.9) circle (5.43130215351701pt);
\fill (0.866025403784439, -0.5) circle (5.43130215351701pt);
\fill (0, 1) circle (5.43130215351701pt);
\fill (-0.866025403784438, -0.5) circle (5.43130215351701pt);
\fill (-2.51147367097487, -1.45) circle (5.43130215351701pt);
\fill (-1.73205080756888, -1) circle (5.43130215351701pt);
\fill (0, 2) circle (5.43130215351701pt);
\end{tikzpicture}

%% file: 1163550515.tikz
\begin{tikzpicture}[scale=0.4]
\draw (2,0) -- (1,0);
\draw (2,0) -- (2.9,0);
\draw (1,0) -- (-0.5,0.866025403784439);
\draw (1,0) -- (-0.5,-0.866025403784438);
\draw (-0.5,0.866025403784439) -- (-1,1.73205080756888);
\draw (-0.5,0.866025403784439) -- (-0.5,-0.866025403784438);
\draw (-1,1.73205080756888) -- (-1.45,2.51147367097487);
\draw (2.9,0) -- (3.71,0);
\draw (-1.45,2.51147367097487) -- (-1.855,3.21295424804027);
\fill (3.71, 0) circle (5.43130215351701pt);
\fill (-0.5, -0.866025403784438) circle (5.43130215351701pt);
\fill (2, 0) circle (5.43130215351701pt);
\fill (-1, 1.73205080756888) circle (5.43130215351701pt);
\fill (-0.5, 0.866025403784439) circle (5.43130215351701pt);
\fill (2.9, 0) circle (5.43130215351701pt);
\fill (-1.855, 3.21295424804027) circle (5.43130215351701pt);
\fill (-1.45, 2.51147367097487) circle (5.43130215351701pt);
\fill (1, 0) circle (5.43130215351701pt);
\end{tikzpicture}

%% file: 1702484648.tikz
\begin{tikzpicture}[scale=0.4]
\draw (2,0) -- (1,0);
\draw (2,0) -- (2.9,0);
\draw (1,0) -- (6.12323399573677e-17,1);
\draw (1,0) -- (-1.83697019872103e-16,-1);
\draw (6.12323399573677e-17,1) -- (-1,1.22464679914735e-16);
\draw (-1,1.22464679914735e-16) -- (-1.83697019872103e-16,-1);
\draw (2.9,0) -- (3.71,0);
\draw (-1.83697019872103e-16,-1) -- (-3.67394039744206e-16,-2);
\draw (-3.67394039744206e-16,-2) -- (-5.32721357629099e-16,-2.9);
\fill (6.12323399573677e-17, 1) circle (4.75pt);
\fill (-1, 1.22464679914735e-16) circle (4.75pt);
\fill (2, 0) circle (4.75pt);
\fill (-1.83697019872103e-16, -1) circle (4.75pt);
\fill (-3.67394039744206e-16, -2) circle (4.75pt);
\fill (1, 0) circle (4.75pt);
\fill (3.71, 0) circle (4.75pt);
\fill (-5.32721357629099e-16, -2.9) circle (4.75pt);
\fill (2.9, 0) circle (4.75pt);
\end{tikzpicture}

%% file: 24790958.tikz
\begin{tikzpicture}[scale=0.4]
\draw (2,0) -- (1,0);
\draw (2,0) -- (2.9,0);
\draw (1,0) -- (6.12323399573677e-17,1);
\draw (1,0) -- (-1.83697019872103e-16,-1);
\draw (6.12323399573677e-17,1) -- (-1,1.22464679914735e-16);
\draw (-1,1.22464679914735e-16) -- (-1.83697019872103e-16,-1);
\draw (-1,1.22464679914735e-16) -- (-2,2.44929359829471e-16);
\draw (2.9,0) -- (3.71,0);
\draw (-2,2.44929359829471e-16) -- (-2.9,3.55147571752732e-16);
\fill (-2, 2.44929359829471e-16) circle (4.75pt);
\fill (1, 0) circle (4.75pt);
\fill (-2.9, 3.55147571752732e-16) circle (4.75pt);
\fill (2.9, 0) circle (4.75pt);
\fill (6.12323399573677e-17, 1) circle (4.75pt);
\fill (-1, 1.22464679914735e-16) circle (4.75pt);
\fill (2, 0) circle (4.75pt);
\fill (-1.83697019872103e-16, -1) circle (4.75pt);
\fill (3.71, 0) circle (4.75pt);
\end{tikzpicture}

%% file: 293328903.tikz
\begin{tikzpicture}[scale=0.4]
\draw (0,2.9) -- (0,2);
\draw (0,2) -- (0,1);
\draw (0,1) -- (1,6.12323399573677e-17);
\draw (0,1) -- (-1,-1.83697019872103e-16);
\draw (1,6.12323399573677e-17) -- (1.22464679914735e-16,-1);
\draw (1,6.12323399573677e-17) -- (2,1.22464679914735e-16);
\draw (-1,-1.83697019872103e-16) -- (1.22464679914735e-16,-1);
\draw (-1,-1.83697019872103e-16) -- (-2,-3.67394039744206e-16);
\draw (-2,-3.67394039744206e-16) -- (-2.9,-5.32721357629099e-16);
\fill (0, 2.9) circle (4.75pt);
\fill (1.22464679914735e-16, -1) circle (4.75pt);
\fill (2, 1.22464679914735e-16) circle (4.75pt);
\fill (0, 1) circle (4.75pt);
\fill (1, 6.12323399573677e-17) circle (4.75pt);
\fill (-1, -1.83697019872103e-16) circle (4.75pt);
\fill (0, 2) circle (4.75pt);
\fill (-2, -3.67394039744206e-16) circle (4.75pt);
\fill (-2.9, -5.32721357629099e-16) circle (4.75pt);
\end{tikzpicture}

%% file: 27921636.tikz
\begin{tikzpicture}[scale=0.4]
\draw (0,1) -- (0.951056516295154,0.309016994374947);
\draw (0,1) -- (-0.951056516295154,0.309016994374947);
\draw (0.951056516295154,0.309016994374947) -- (0.587785252292473,-0.809016994374947);
\draw (-0.951056516295154,0.309016994374947) -- (-1.90211303259031,0.618033988749894);
\draw (-0.951056516295154,0.309016994374947) -- (-0.587785252292473,-0.809016994374947);
\draw (0.587785252292473,-0.809016994374947) -- (1.17557050458495,-1.61803398874989);
\draw (0.587785252292473,-0.809016994374947) -- (-0.587785252292473,-0.809016994374947);
\draw (-1.90211303259031,0.618033988749894) -- (-2.75806389725595,0.896149283687347);
\draw (1.17557050458495,-1.61803398874989) -- (1.70457723164817,-2.34614928368735);
\fill (0, 1) circle (4.26520876399966pt);
\fill (1.17557050458495, -1.61803398874989) circle (4.26520876399966pt);
\fill (-2.75806389725595, 0.896149283687347) circle (4.26520876399966pt);
\fill (-0.951056516295154, 0.309016994374947) circle (4.26520876399966pt);
\fill (0.587785252292473, -0.809016994374947) circle (4.26520876399966pt);
\fill (-1.90211303259031, 0.618033988749894) circle (4.26520876399966pt);
\fill (0.951056516295154, 0.309016994374947) circle (4.26520876399966pt);
\fill (-0.587785252292473, -0.809016994374947) circle (4.26520876399966pt);
\fill (1.70457723164817, -2.34614928368735) circle (4.26520876399966pt);
\end{tikzpicture}

%% file: 552447728.tikz
\begin{tikzpicture}[scale=0.4]
\draw (0,2) -- (0,1);
\draw (0,1) -- (0.951056516295154,0.309016994374947);
\draw (0,1) -- (-0.951056516295154,0.309016994374947);
\draw (0.951056516295154,0.309016994374947) -- (0.587785252292473,-0.809016994374947);
\draw (-0.951056516295154,0.309016994374947) -- (-0.587785252292473,-0.809016994374947);
\draw (-0.951056516295154,0.309016994374947) -- (-1.90211303259031,0.618033988749894);
\draw (0.587785252292473,-0.809016994374947) -- (-0.587785252292473,-0.809016994374947);
\draw (0.587785252292473,-0.809016994374947) -- (1.17557050458495,-1.61803398874989);
\draw (-1.90211303259031,0.618033988749894) -- (-2.75806389725595,0.896149283687347);
\fill (-0.951056516295154, 0.309016994374947) circle (4.26520876399966pt);
\fill (0.951056516295154, 0.309016994374947) circle (4.26520876399966pt);
\fill (1.17557050458495, -1.61803398874989) circle (4.26520876399966pt);
\fill (0, 2) circle (4.26520876399966pt);
\fill (-0.587785252292473, -0.809016994374947) circle (4.26520876399966pt);
\fill (-2.75806389725595, 0.896149283687347) circle (4.26520876399966pt);
\fill (-1.90211303259031, 0.618033988749894) circle (4.26520876399966pt);
\fill (0, 1) circle (4.26520876399966pt);
\fill (0.587785252292473, -0.809016994374947) circle (4.26520876399966pt);
\end{tikzpicture}

%% file: 1837752289.tikz
\begin{tikzpicture}[scale=0.4]
\draw (0,2.9) -- (0,2);
\draw (0,2) -- (0,1);
\draw (0,1) -- (0.951056516295154,0.309016994374947);
\draw (0,1) -- (-0.951056516295154,0.309016994374947);
\draw (0.951056516295154,0.309016994374947) -- (0.587785252292473,-0.809016994374947);
\draw (0.587785252292473,-0.809016994374947) -- (-0.587785252292473,-0.809016994374947);
\draw (-0.951056516295154,0.309016994374947) -- (-0.587785252292473,-0.809016994374947);
\draw (-0.951056516295154,0.309016994374947) -- (-1.90211303259031,0.618033988749894);
\draw (-1.90211303259031,0.618033988749894) -- (-2.75806389725595,0.896149283687347);
\fill (0, 2.9) circle (4.26520876399966pt);
\fill (-0.951056516295154, 0.309016994374947) circle (4.26520876399966pt);
\fill (-1.90211303259031, 0.618033988749894) circle (4.26520876399966pt);
\fill (0, 1) circle (4.26520876399966pt);
\fill (0.951056516295154, 0.309016994374947) circle (4.26520876399966pt);
\fill (0.587785252292473, -0.809016994374947) circle (4.26520876399966pt);
\fill (0, 2) circle (4.26520876399966pt);
\fill (-0.587785252292473, -0.809016994374947) circle (4.26520876399966pt);
\fill (-2.75806389725595, 0.896149283687347) circle (4.26520876399966pt);
\end{tikzpicture}

%% file: 562584937.tikz
\begin{tikzpicture}[scale=0.4]
\draw (2,0) -- (1,0);
\draw (2,0) -- (2.9,0);
\draw (1,0) -- (0.309016994374947,0.951056516295154);
\draw (1,0) -- (0.309016994374947,-0.951056516295154);
\draw (0.309016994374947,0.951056516295154) -- (-0.809016994374947,0.587785252292473);
\draw (-0.809016994374947,0.587785252292473) -- (-0.809016994374947,-0.587785252292473);
\draw (2.9,0) -- (3.71,0);
\draw (-0.809016994374947,-0.587785252292473) -- (0.309016994374947,-0.951056516295154);
\draw (-0.809016994374947,-0.587785252292473) -- (-1.61803398874989,-1.17557050458495);
\fill (3.71, 0) circle (4.26520876399966pt);
\fill (1, 0) circle (4.26520876399966pt);
\fill (-1.61803398874989, -1.17557050458495) circle (4.26520876399966pt);
\fill (2.9, 0) circle (4.26520876399966pt);
\fill (0.309016994374947, 0.951056516295154) circle (4.26520876399966pt);
\fill (-0.809016994374947, 0.587785252292473) circle (4.26520876399966pt);
\fill (2, 0) circle (4.26520876399966pt);
\fill (-0.809016994374947, -0.587785252292473) circle (4.26520876399966pt);
\fill (0.309016994374947, -0.951056516295154) circle (4.26520876399966pt);
\end{tikzpicture}

%% file: 827130457.tikz
\begin{tikzpicture}[scale=0.4]
\draw (0,1) -- (0.866025403784439,0.5);
\draw (0,1) -- (-0.866025403784439,0.5);
\draw (0.866025403784439,0.5) -- (1.73205080756888,1);
\draw (0.866025403784439,0.5) -- (0.866025403784439,-0.5);
\draw (1.73205080756888,1) -- (2.51147367097487,1.45);
\draw (-0.866025403784439,0.5) -- (-0.866025403784438,-0.5);
\draw (2.51147367097487,1.45) -- (3.21295424804027,1.855);
\draw (-0.866025403784438,-0.5) -- (1.22464679914735e-16,-1);
\draw (1.22464679914735e-16,-1) -- (0.866025403784439,-0.5);
\fill (1.22464679914735e-16, -1) circle (3.8909863237109pt);
\fill (-0.866025403784438, -0.5) circle (3.8909863237109pt);
\fill (0, 1) circle (3.8909863237109pt);
\fill (-0.866025403784439, 0.5) circle (3.8909863237109pt);
\fill (1.73205080756888, 1) circle (3.8909863237109pt);
\fill (2.51147367097487, 1.45) circle (3.8909863237109pt);
\fill (0.866025403784439, -0.5) circle (3.8909863237109pt);
\fill (0.866025403784439, 0.5) circle (3.8909863237109pt);
\fill (3.21295424804027, 1.855) circle (3.8909863237109pt);
\end{tikzpicture}

%% file: 877228084.tikz
\begin{tikzpicture}[scale=0.4]
\draw (0,2) -- (0,1);
\draw (0,1) -- (0.866025403784439,0.5);
\draw (0,1) -- (-0.866025403784439,0.5);
\draw (0.866025403784439,0.5) -- (0.866025403784439,-0.5);
\draw (0.866025403784439,-0.5) -- (1.22464679914735e-16,-1);
\draw (-0.866025403784439,0.5) -- (-0.866025403784438,-0.5);
\draw (-0.866025403784439,0.5) -- (-1.73205080756888,1);
\draw (1.22464679914735e-16,-1) -- (-0.866025403784438,-0.5);
\draw (-1.73205080756888,1) -- (-2.51147367097487,1.45);
\fill (-2.51147367097487, 1.45) circle (3.8909863237109pt);
\fill (-0.866025403784438, -0.5) circle (3.8909863237109pt);
\fill (0, 1) circle (3.8909863237109pt);
\fill (-0.866025403784439, 0.5) circle (3.8909863237109pt);
\fill (0.866025403784439, -0.5) circle (3.8909863237109pt);
\fill (0.866025403784439, 0.5) circle (3.8909863237109pt);
\fill (-1.73205080756888, 1) circle (3.8909863237109pt);
\fill (1.22464679914735e-16, -1) circle (3.8909863237109pt);
\fill (0, 2) circle (3.8909863237109pt);
\end{tikzpicture}

%% file: 1095260065.tikz
\begin{tikzpicture}[scale=0.4]
\draw (0,2) -- (0,1);
\draw (0,1) -- (0.866025403784439,0.5);
\draw (0,1) -- (-0.866025403784439,0.5);
\draw (0.866025403784439,0.5) -- (0.866025403784439,-0.5);
\draw (-0.866025403784439,0.5) -- (-0.866025403784438,-0.5);
\draw (1.22464679914735e-16,-1) -- (0.866025403784439,-0.5);
\draw (1.22464679914735e-16,-1) -- (-0.866025403784438,-0.5);
\draw (0.866025403784439,-0.5) -- (1.73205080756888,-1);
\draw (-0.866025403784438,-0.5) -- (-1.73205080756888,-1);
\fill (1.22464679914735e-16, -1) circle (3.8909863237109pt);
\fill (-1.73205080756888, -1) circle (3.8909863237109pt);
\fill (0, 1) circle (3.8909863237109pt);
\fill (-0.866025403784438, -0.5) circle (3.8909863237109pt);
\fill (1.73205080756888, -1) circle (3.8909863237109pt);
\fill (0, 2) circle (3.8909863237109pt);
\fill (0.866025403784439, -0.5) circle (3.8909863237109pt);
\fill (-0.866025403784439, 0.5) circle (3.8909863237109pt);
\fill (0.866025403784439, 0.5) circle (3.8909863237109pt);
\end{tikzpicture}

%% file: 1347018042.tikz
\begin{tikzpicture}[scale=0.4]
\draw (0,1) -- (0.866025403784439,0.5);
\draw (0,1) -- (-0.866025403784439,0.5);
\draw (0.866025403784439,0.5) -- (0.866025403784439,-0.5);
\draw (-0.866025403784439,0.5) -- (-0.866025403784438,-0.5);
\draw (0.866025403784439,-0.5) -- (1.73205080756888,-1);
\draw (0.866025403784439,-0.5) -- (1.22464679914735e-16,-1);
\draw (1.73205080756888,-1) -- (2.51147367097487,-1.45);
\draw (-0.866025403784438,-0.5) -- (1.22464679914735e-16,-1);
\draw (-0.866025403784438,-0.5) -- (-1.73205080756888,-1);
\fill (0.866025403784439, -0.5) circle (3.8909863237109pt);
\fill (-0.866025403784439, 0.5) circle (3.8909863237109pt);
\fill (1.22464679914735e-16, -1) circle (3.8909863237109pt);
\fill (0, 1) circle (3.8909863237109pt);
\fill (-0.866025403784438, -0.5) circle (3.8909863237109pt);
\fill (-1.73205080756888, -1) circle (3.8909863237109pt);
\fill (0.866025403784439, 0.5) circle (3.8909863237109pt);
\fill (2.51147367097487, -1.45) circle (3.8909863237109pt);
\fill (1.73205080756888, -1) circle (3.8909863237109pt);
\end{tikzpicture}

%% file: 2128804136.tikz
\begin{tikzpicture}[scale=0.5]
\draw (0,1) -- (0.78183148246803,0.623489801858734);
\draw (0,1) -- (-0.78183148246803,0.623489801858733);
\draw (0.78183148246803,0.623489801858734) -- (0.974927912181824,-0.222520933956314);
\draw (-0.78183148246803,0.623489801858733) -- (-0.974927912181824,-0.222520933956315);
\draw (0.974927912181824,-0.222520933956314) -- (1.94985582436365,-0.445041867912629);
\draw (0.974927912181824,-0.222520933956314) -- (0.433883739117558,-0.900968867902419);
\draw (1.94985582436365,-0.445041867912629) -- (2.82729094532729,-0.645310708473312);
\draw (-0.974927912181824,-0.222520933956315) -- (-0.433883739117558,-0.900968867902419);
\draw (-0.433883739117558,-0.900968867902419) -- (0.433883739117558,-0.900968867902419);
\fill (-0.433883739117558, -0.900968867902419) circle (3.58621578407382pt);
\fill (-0.974927912181824, -0.222520933956315) circle (3.58621578407382pt);
\fill (0, 1) circle (3.58621578407382pt);
\fill (0.974927912181824, -0.222520933956314) circle (3.58621578407382pt);
\fill (-0.78183148246803, 0.623489801858733) circle (3.58621578407382pt);
\fill (1.94985582436365, -0.445041867912629) circle (3.58621578407382pt);
\fill (0.433883739117558, -0.900968867902419) circle (3.58621578407382pt);
\fill (2.82729094532729, -0.645310708473312) circle (3.58621578407382pt);
\fill (0.78183148246803, 0.623489801858734) circle (3.58621578407382pt);
\end{tikzpicture}

%% file: 967574001.tikz
\begin{tikzpicture}[scale=0.5]
\draw (0,1) -- (0.78183148246803,0.623489801858734);
\draw (0,1) -- (-0.78183148246803,0.623489801858733);
\draw (0.78183148246803,0.623489801858734) -- (0.974927912181824,-0.222520933956314);
\draw (-0.78183148246803,0.623489801858733) -- (-0.974927912181824,-0.222520933956315);
\draw (0.974927912181824,-0.222520933956314) -- (0.433883739117558,-0.900968867902419);
\draw (-0.974927912181824,-0.222520933956315) -- (-0.433883739117558,-0.900968867902419);
\draw (-0.974927912181824,-0.222520933956315) -- (-1.94985582436365,-0.445041867912629);
\draw (0.433883739117558,-0.900968867902419) -- (-0.433883739117558,-0.900968867902419);
\draw (0.433883739117558,-0.900968867902419) -- (0.867767478235116,-1.80193773580484);
\fill (-0.974927912181824, -0.222520933956315) circle (3.58621578407382pt);
\fill (-0.433883739117558, -0.900968867902419) circle (3.58621578407382pt);
\fill (0.78183148246803, 0.623489801858734) circle (3.58621578407382pt);
\fill (0.867767478235116, -1.80193773580484) circle (3.58621578407382pt);
\fill (0, 1) circle (3.58621578407382pt);
\fill (0.433883739117558, -0.900968867902419) circle (3.58621578407382pt);
\fill (-1.94985582436365, -0.445041867912629) circle (3.58621578407382pt);
\fill (-0.78183148246803, 0.623489801858733) circle (3.58621578407382pt);
\fill (0.974927912181824, -0.222520933956314) circle (3.58621578407382pt);
\end{tikzpicture}

%% file: 1623687162.tikz
\begin{tikzpicture}[scale=0.6]
\draw (0,1) -- (0.707106781186547,0.707106781186548);
\draw (0,1) -- (-0.707106781186548,0.707106781186547);
\draw (0.707106781186547,0.707106781186548) -- (1,6.12323399573677e-17);
\draw (-0.707106781186548,0.707106781186547) -- (-1,-1.83697019872103e-16);
\draw (1,6.12323399573677e-17) -- (0.707106781186548,-0.707106781186547);
\draw (-1,-1.83697019872103e-16) -- (-0.707106781186547,-0.707106781186548);
\draw (0.707106781186548,-0.707106781186547) -- (1.22464679914735e-16,-1);
\draw (0.707106781186548,-0.707106781186547) -- (1.4142135623731,-1.41421356237309);
\draw (1.22464679914735e-16,-1) -- (-0.707106781186547,-0.707106781186548);
\fill (-0.707106781186547, -0.707106781186548) circle (3.32842712474619pt);
\fill (0.707106781186547, 0.707106781186548) circle (3.32842712474619pt);
\fill (1.22464679914735e-16, -1) circle (3.32842712474619pt);
\fill (-1, -1.83697019872103e-16) circle (3.32842712474619pt);
\fill (1, 6.12323399573677e-17) circle (3.32842712474619pt);
\fill (-0.707106781186548, 0.707106781186547) circle (3.32842712474619pt);
\fill (1.4142135623731, -1.41421356237309) circle (3.32842712474619pt);
\fill (0.707106781186548, -0.707106781186547) circle (3.32842712474619pt);
\fill (0, 1) circle (3.32842712474619pt);
\end{tikzpicture}

%% file: 1341331914.tikz
\begin{tikzpicture}[scale=0.6]
\draw (0,1) -- (0.642787609686539,0.766044443118978);
\draw (0,1) -- (-0.64278760968654,0.766044443118978);
\draw (0.642787609686539,0.766044443118978) -- (0.984807753012208,0.17364817766693);
\draw (-0.64278760968654,0.766044443118978) -- (-0.984807753012208,0.17364817766693);
\draw (0.984807753012208,0.17364817766693) -- (0.866025403784439,-0.5);
\draw (-0.984807753012208,0.17364817766693) -- (-0.866025403784438,-0.5);
\draw (0.866025403784439,-0.5) -- (0.342020143325669,-0.939692620785908);
\draw (-0.342020143325669,-0.939692620785908) -- (-0.866025403784438,-0.5);
\draw (-0.342020143325669,-0.939692620785908) -- (0.342020143325669,-0.939692620785908);
\fill (0.342020143325669, -0.939692620785908) circle (3.10416666666667pt);
\fill (0.642787609686539, 0.766044443118978) circle (3.10416666666667pt);
\fill (-0.342020143325669, -0.939692620785908) circle (3.10416666666667pt);
\fill (-0.984807753012208, 0.17364817766693) circle (3.10416666666667pt);
\fill (0.984807753012208, 0.17364817766693) circle (3.10416666666667pt);
\fill (-0.64278760968654, 0.766044443118978) circle (3.10416666666667pt);
\fill (-0.866025403784438, -0.5) circle (3.10416666666667pt);
\fill (0, 1) circle (3.10416666666667pt);
\fill (0.866025403784439, -0.5) circle (3.10416666666667pt);
\end{tikzpicture}

%% file: minkowski3.tex
\begin{tikzpicture}[x  = {(0.9cm,-0.076cm)},
                    y  = {(-0.06cm,0.95cm)},
                    z  = {(-0.44cm,-0.29cm)},
                    scale = 1,
                    color = {lightgray}]


  \coordinate (v0__1) at (0, 0, 0);
  \coordinate (v1__1) at (0, 0, 1);
  \coordinate (v2__1) at (0, 1, 0);
  \coordinate (v3__1) at (1, 0, 0);

  \definecolor{vertexcolor__1}{rgb}{ 0 0 0 }

  \tikzstyle{vertexstyle__1_0} = [circle, scale=0.38, fill=vertexcolor__1,]
  \tikzstyle{vertexstyle__1_1} = [circle, scale=0.38, fill=vertexcolor__1,]
  \tikzstyle{vertexstyle__1_2} = [circle, scale=0.38, fill=vertexcolor__1,]
  \tikzstyle{vertexstyle__1_3} = [circle, scale=0.38, fill=vertexcolor__1,]

  \definecolor{facetcolor__1}{rgb}{ 1 1 1 }

  \definecolor{edgecolor__1}{rgb}{ 0 0 0 }
  \tikzstyle{facetstyle__1} = [fill=facetcolor__1, fill opacity=0.25, draw=edgecolor__1, line width=1 pt, line cap=round, line join=round]

  \draw[facetstyle__1] (v0__1) -- (v1__1) -- (v2__1) -- (v0__1) -- cycle;
  \draw[facetstyle__1] (v3__1) -- (v1__1) -- (v0__1) -- (v3__1) -- cycle;
  \draw[facetstyle__1] (v0__1) -- (v2__1) -- (v3__1) -- (v0__1) -- cycle;

   \node at (v0__1) [vertexstyle__1_0] {};

  \draw[facetstyle__1] (v2__1) -- (v1__1) -- (v3__1) -- (v2__1) -- cycle;

  \foreach \i in {1,2,3} {
    \node at (v\i__1) [vertexstyle__1_\i] {};
  }

  \coordinate (v0__2) at (0, 0, 1);
  \coordinate (v1__2) at (0, 0, 2);
  \coordinate (v2__2) at (0, 1, 0);
  \coordinate (v3__2) at (0, 1, 1);
  \coordinate (v4__2) at (1, 0, 0);
  \coordinate (v5__2) at (1, 0, 1);

  \definecolor{vertexcolor__2}{rgb}{ 0 0 0 }

  \tikzstyle{vertexstyle__2_0} = [circle, scale=0.38, fill=vertexcolor__2,]
  \tikzstyle{vertexstyle__2_1} = [circle, scale=0.38, fill=vertexcolor__2,]
  \tikzstyle{vertexstyle__2_2} = [circle, scale=0.38, fill=vertexcolor__2,]
  \tikzstyle{vertexstyle__2_3} = [circle, scale=0.38, fill=vertexcolor__2,]
  \tikzstyle{vertexstyle__2_4} = [circle, scale=0.38, fill=vertexcolor__2,]
  \tikzstyle{vertexstyle__2_5} = [circle, scale=0.38, fill=vertexcolor__2,]

  \definecolor{facetcolor__2}{rgb}{ 0.57843137254 0.74705882352 0.90196078431 }

  \definecolor{edgecolor__2}{rgb}{ 0 0 0 }
  \tikzstyle{facetstyle__2} = [fill=facetcolor__2, fill opacity=0.6, draw=edgecolor__2, line width=1 pt, line cap=round, line join=round]

  \draw[facetstyle__2] (v4__2) -- (v5__2) -- (v1__2) -- (v0__2) -- (v4__2) -- cycle;
  \draw[facetstyle__2] (v4__2) -- (v0__2) -- (v2__2) -- (v4__2) -- cycle;
  \draw[facetstyle__2] (v0__2) -- (v1__2) -- (v3__2) -- (v2__2) -- (v0__2) -- cycle;

   \node at (v0__2) [vertexstyle__2_0] {};

  \draw[facetstyle__2] (v3__2) -- (v5__2) -- (v4__2) -- (v2__2) -- (v3__2) -- cycle;
  \draw[facetstyle__2] (v1__2) -- (v5__2) -- (v3__2) -- (v1__2) -- cycle;

  \foreach \i in {1,3,5,2,4} {
    \node at (v\i__2) [vertexstyle__2_\i] {};
  }

  \coordinate (v0__8) at (0, 1, 0);
  \coordinate (v1__8) at (0, 1, 1);
  \coordinate (v2__8) at (0, 2, 0);
  \coordinate (v3__8) at (1, 0, 0);
  \coordinate (v4__8) at (1, 0, 1);
  \coordinate (v5__8) at (1, 1, 0);

  \definecolor{vertexcolor__8}{rgb}{ 0 0 0 }

  \tikzstyle{vertexstyle__8_0} = [circle, scale=0.38, fill=vertexcolor__8,]
  \tikzstyle{vertexstyle__8_1} = [circle, scale=0.38, fill=vertexcolor__8,]
  \tikzstyle{vertexstyle__8_2} = [circle, scale=0.38, fill=vertexcolor__8,]
  \tikzstyle{vertexstyle__8_3} = [circle, scale=0.38, fill=vertexcolor__8,]
  \tikzstyle{vertexstyle__8_4} = [circle, scale=0.38, fill=vertexcolor__8,]
  \tikzstyle{vertexstyle__8_5} = [circle, scale=0.38, fill=vertexcolor__8,]

  \definecolor{facetcolor__8}{rgb}{ 0.67843137254 0.84705882352 0.90196078431 }

  \definecolor{edgecolor__8}{rgb}{ 0 0 0 }
  \tikzstyle{facetstyle__8} = [fill=facetcolor__8, fill opacity=0.6, draw=edgecolor__8, line width=1 pt, line cap=round, line join=round]

  \draw[facetstyle__8] (v2__8) -- (v5__8) -- (v3__8) -- (v0__8) -- (v2__8) -- cycle;
  \draw[facetstyle__8] (v0__8) -- (v3__8) -- (v4__8) -- (v1__8) -- (v0__8) -- cycle;
  \draw[facetstyle__8] (v2__8) -- (v0__8) -- (v1__8) -- (v2__8) -- cycle;

   \node at (v0__8) [vertexstyle__8_0] {};

  \draw[facetstyle__8] (v4__8) -- (v5__8) -- (v2__8) -- (v1__8) -- (v4__8) -- cycle;
  \draw[facetstyle__8] (v3__8) -- (v5__8) -- (v4__8) -- (v3__8) -- cycle;

  \foreach \i in {1,4,3,2,5} {
    \node at (v\i__8) [vertexstyle__8_\i] {};
  }

  \coordinate (v0__10) at (0, 1, 1);
  \coordinate (v1__10) at (0, 2, 0);
  \coordinate (v2__10) at (1, 0, 1);
  \coordinate (v3__10) at (1, 1, 0);
  \coordinate (v4__10) at (1, 1, 1);
  \coordinate (v5__10) at (1, 2, 0);

  \definecolor{vertexcolor__10}{rgb}{ 0 0 0 }

  \tikzstyle{vertexstyle__10_0} = [circle, scale=0.38, fill=vertexcolor__10,]
  \tikzstyle{vertexstyle__10_1} = [circle, scale=0.38, fill=vertexcolor__10,]
  \tikzstyle{vertexstyle__10_2} = [circle, scale=0.38, fill=vertexcolor__10,]
  \tikzstyle{vertexstyle__10_3} = [circle, scale=0.38, fill=vertexcolor__10,]
  \tikzstyle{vertexstyle__10_4} = [circle, scale=0.38, fill=vertexcolor__10,]
  \tikzstyle{vertexstyle__10_5} = [circle, scale=0.38, fill=vertexcolor__10,]

  \definecolor{facetcolor__10}{rgb}{ 0.67843137254 0.84705882352 0.90196078431 }

  \definecolor{edgecolor__10}{rgb}{ 0 0 0 }
  \tikzstyle{facetstyle__10} = [fill=facetcolor__10, fill opacity=0.6, draw=edgecolor__10, line width=1 pt, line cap=round, line join=round]

  \draw[facetstyle__10] (v0__10) -- (v1__10) -- (v3__10) -- (v2__10) -- (v0__10) -- cycle;
  \draw[facetstyle__10] (v1__10) -- (v5__10) -- (v3__10) -- (v1__10) -- cycle;
  \draw[facetstyle__10] (v4__10) -- (v5__10) -- (v1__10) -- (v0__10) -- (v4__10) -- cycle;
  \draw[facetstyle__10] (v3__10) -- (v5__10) -- (v4__10) -- (v2__10) -- (v3__10) -- cycle;
  \draw[facetstyle__10] (v4__10) -- (v0__10) -- (v2__10) -- (v4__10) -- cycle;

  \foreach \i in {0,2,4,1,3,5} {
    \node at (v\i__10) [vertexstyle__10_\i] {};
  }

  \coordinate (v0__6) at (0, 0, 2);
  \coordinate (v1__6) at (0, 1, 1);
  \coordinate (v2__6) at (1, 0, 1);
  \coordinate (v3__6) at (1, 1, 1);

  \definecolor{vertexcolor__6}{rgb}{ 0 0 0 }

  \tikzstyle{vertexstyle__6_0} = [circle, scale=0.38, fill=vertexcolor__6,]
  \tikzstyle{vertexstyle__6_1} = [circle, scale=0.38, fill=vertexcolor__6,]
  \tikzstyle{vertexstyle__6_2} = [circle, scale=0.38, fill=vertexcolor__6,]
  \tikzstyle{vertexstyle__6_3} = [circle, scale=0.38, fill=vertexcolor__6,]

  \definecolor{facetcolor__6}{rgb}{ 1 1 1 }

  \definecolor{edgecolor__6}{rgb}{ 0 0 0 }
  \tikzstyle{facetstyle__6} = [fill=facetcolor__6, fill opacity=0.25, draw=edgecolor__6, line width=1 pt, line cap=round, line join=round]

  \draw[facetstyle__6] (v0__6) -- (v1__6) -- (v2__6) -- (v0__6) -- cycle;
  \draw[facetstyle__6] (v2__6) -- (v1__6) -- (v3__6) -- (v2__6) -- cycle;
  \draw[facetstyle__6] (v3__6) -- (v1__6) -- (v0__6) -- (v3__6) -- cycle;
  \draw[facetstyle__6] (v0__6) -- (v2__6) -- (v3__6) -- (v0__6) -- cycle;

  \foreach \i in {0,1,2,3} {
    \node at (v\i__6) [vertexstyle__6_\i] {};
  }

  \coordinate (v0__5) at (0, 0, 2);
  \coordinate (v1__5) at (0, 1, 1);
  \coordinate (v2__5) at (0, 1, 2);
  \coordinate (v3__5) at (1, 1, 1);

  \definecolor{vertexcolor__5}{rgb}{ 0 0 0 }

  \tikzstyle{vertexstyle__5_0} = [circle, scale=0.38, fill=vertexcolor__5,]
  \tikzstyle{vertexstyle__5_1} = [circle, scale=0.38, fill=vertexcolor__5,]
  \tikzstyle{vertexstyle__5_2} = [circle, scale=0.38, fill=vertexcolor__5,]
  \tikzstyle{vertexstyle__5_3} = [circle, scale=0.38, fill=vertexcolor__5,]

  \definecolor{facetcolor__5}{rgb}{ 1 1 1 }

  \definecolor{edgecolor__5}{rgb}{ 0 0 0 }
  \tikzstyle{facetstyle__5} = [fill=facetcolor__5, fill opacity=0.25, draw=edgecolor__5, line width=1 pt, line cap=round, line join=round]

  \draw[facetstyle__5] (v1__5) -- (v0__5) -- (v2__5) -- (v1__5) -- cycle;
  \draw[facetstyle__5] (v3__5) -- (v0__5) -- (v1__5) -- (v3__5) -- cycle;
  \draw[facetstyle__5] (v2__5) -- (v0__5) -- (v3__5) -- (v2__5) -- cycle;
  \draw[facetstyle__5] (v1__5) -- (v2__5) -- (v3__5) -- (v1__5) -- cycle;

  \foreach \i in {0,2,1,3} {
    \node at (v\i__5) [vertexstyle__5_\i] {};
  }

  \coordinate (v0__7) at (0, 0, 2);
  \coordinate (v1__7) at (1, 0, 1);
  \coordinate (v2__7) at (1, 0, 2);
  \coordinate (v3__7) at (1, 1, 1);

  \definecolor{vertexcolor__7}{rgb}{ 0 0 0 }

  \tikzstyle{vertexstyle__7_0} = [circle, scale=0.38, fill=vertexcolor__7,]
  \tikzstyle{vertexstyle__7_1} = [circle, scale=0.38, fill=vertexcolor__7,]
  \tikzstyle{vertexstyle__7_2} = [circle, scale=0.38, fill=vertexcolor__7,]
  \tikzstyle{vertexstyle__7_3} = [circle, scale=0.38, fill=vertexcolor__7,]

  \definecolor{facetcolor__7}{rgb}{ 1 1 1 }

  \definecolor{edgecolor__7}{rgb}{ 0 0 0 }
  \tikzstyle{facetstyle__7} = [fill=facetcolor__7, fill opacity=0.25, draw=edgecolor__7, line width=1 pt, line cap=round, line join=round]

  \draw[facetstyle__7] (v0__7) -- (v1__7) -- (v2__7) -- (v0__7) -- cycle;
  \draw[facetstyle__7] (v3__7) -- (v1__7) -- (v0__7) -- (v3__7) -- cycle;
  \draw[facetstyle__7] (v0__7) -- (v2__7) -- (v3__7) -- (v0__7) -- cycle;
  \draw[facetstyle__7] (v2__7) -- (v1__7) -- (v3__7) -- (v2__7) -- cycle;

   \node at (v0__7) [vertexstyle__7_0] {};

  \draw[facetstyle__7] (v2__7) -- (v1__7) -- (v3__7) -- (v2__7) -- cycle;

  \foreach \i in {0,2,1,3} {
    \node at (v\i__7) [vertexstyle__7_\i] {};
  }

  \coordinate (v0__12) at (1, 0, 0);
  \coordinate (v1__12) at (1, 0, 1);
  \coordinate (v2__12) at (1, 1, 0);
  \coordinate (v3__12) at (2, 0, 0);

  \definecolor{vertexcolor__12}{rgb}{ 0 0 0 }

  \tikzstyle{vertexstyle__12_0} = [circle, scale=0.38, fill=vertexcolor__12,]
  \tikzstyle{vertexstyle__12_1} = [circle, scale=0.38, fill=vertexcolor__12,]
  \tikzstyle{vertexstyle__12_2} = [circle, scale=0.38, fill=vertexcolor__12,]
  \tikzstyle{vertexstyle__12_3} = [circle, scale=0.38, fill=vertexcolor__12,]

  \definecolor{facetcolor__12}{rgb}{ 1 1 1 }

  \definecolor{edgecolor__12}{rgb}{ 0 0 0 }
  \tikzstyle{facetstyle__12} = [fill=facetcolor__12, fill opacity=0.25, draw=edgecolor__12, line width=1 pt, line cap=round, line join=round]

  \draw[facetstyle__12] (v0__12) -- (v1__12) -- (v2__12) -- (v0__12) -- cycle;
  \draw[facetstyle__12] (v3__12) -- (v1__12) -- (v0__12) -- (v3__12) -- cycle;
  \draw[facetstyle__12] (v0__12) -- (v2__12) -- (v3__12) -- (v0__12) -- cycle;

   \node at (v0__12) [vertexstyle__12_0] {};

  \draw[facetstyle__12] (v2__12) -- (v1__12) -- (v3__12) -- (v2__12) -- cycle;

  \foreach \i in {1,2,3} {
    \node at (v\i__12) [vertexstyle__12_\i] {};
  }

  \coordinate (v0__14) at (1, 0, 1);
  \coordinate (v1__14) at (1, 1, 0);
  \coordinate (v2__14) at (1, 1, 1);
  \coordinate (v3__14) at (1, 2, 0);
  \coordinate (v4__14) at (2, 0, 0);
  \coordinate (v5__14) at (2, 1, 0);

  \definecolor{vertexcolor__14}{rgb}{ 0 0 0 }

  \tikzstyle{vertexstyle__14_0} = [circle, scale=0.38, fill=vertexcolor__14,]
  \tikzstyle{vertexstyle__14_1} = [circle, scale=0.38, fill=vertexcolor__14,]
  \tikzstyle{vertexstyle__14_2} = [circle, scale=0.38, fill=vertexcolor__14,]
  \tikzstyle{vertexstyle__14_3} = [circle, scale=0.38, fill=vertexcolor__14,]
  \tikzstyle{vertexstyle__14_4} = [circle, scale=0.38, fill=vertexcolor__14,]
  \tikzstyle{vertexstyle__14_5} = [circle, scale=0.38, fill=vertexcolor__14,]

  \definecolor{facetcolor__14}{rgb}{ 0.67843137254 0.84705882352 0.90196078431 }

  \definecolor{edgecolor__14}{rgb}{ 0 0 0 }
  \tikzstyle{facetstyle__14} = [fill=facetcolor__14, fill opacity=0.6, draw=edgecolor__14, line width=1 pt, line cap=round, line join=round]

  \draw[facetstyle__14] (v3__14) -- (v5__14) -- (v4__14) -- (v1__14) -- (v3__14) -- cycle;
  \draw[facetstyle__14] (v0__14) -- (v2__14) -- (v3__14) -- (v1__14) -- (v0__14) -- cycle;
  \draw[facetstyle__14] (v4__14) -- (v0__14) -- (v1__14) -- (v4__14) -- cycle;

   \node at (v1__14) [vertexstyle__14_1] {};

  \draw[facetstyle__14] (v4__14) -- (v5__14) -- (v2__14) -- (v0__14) -- (v4__14) -- cycle;
  \draw[facetstyle__14] (v2__14) -- (v5__14) -- (v3__14) -- (v2__14) -- cycle;

  \foreach \i in {0,2,4,3,5} {
    \node at (v\i__14) [vertexstyle__14_\i] {};
  }

  \coordinate (v0__15) at (1, 0, 1);
  \coordinate (v1__15) at (1, 1, 1);
  \coordinate (v2__15) at (2, 0, 0);
  \coordinate (v3__15) at (2, 0, 1);
  \coordinate (v4__15) at (2, 1, 0);
  \coordinate (v5__15) at (3, 0, 0);

  \definecolor{vertexcolor__15}{rgb}{ 0 0 0 }

  \tikzstyle{vertexstyle__15_0} = [circle, scale=0.38, fill=vertexcolor__15,]
  \tikzstyle{vertexstyle__15_1} = [circle, scale=0.38, fill=vertexcolor__15,]
  \tikzstyle{vertexstyle__15_2} = [circle, scale=0.38, fill=vertexcolor__15,]
  \tikzstyle{vertexstyle__15_3} = [circle, scale=0.38, fill=vertexcolor__15,]
  \tikzstyle{vertexstyle__15_4} = [circle, scale=0.38, fill=vertexcolor__15,]
  \tikzstyle{vertexstyle__15_5} = [circle, scale=0.38, fill=vertexcolor__15,]

  \definecolor{facetcolor__15}{rgb}{ 0.67843137254 0.84705882352 0.90196078431 }

  \definecolor{edgecolor__15}{rgb}{ 0 0 0 }
  \tikzstyle{facetstyle__15} = [fill=facetcolor__15, fill opacity=0.6, draw=edgecolor__15, line width=1 pt, line cap=round, line join=round]

  \draw[facetstyle__15] (v5__15) -- (v3__15) -- (v0__15) -- (v2__15) -- (v5__15) -- cycle;
  \draw[facetstyle__15] (v2__15) -- (v0__15) -- (v1__15) -- (v4__15) -- (v2__15) -- cycle;
  \draw[facetstyle__15] (v5__15) -- (v2__15) -- (v4__15) -- (v5__15) -- cycle;

   \node at (v2__15) [vertexstyle__15_2] {};

  \draw[facetstyle__15] (v1__15) -- (v3__15) -- (v5__15) -- (v4__15) -- (v1__15) -- cycle;
  \draw[facetstyle__15] (v0__15) -- (v3__15) -- (v1__15) -- (v0__15) -- cycle;

  \foreach \i in {0,1,3,4,5} {
    \node at (v\i__15) [vertexstyle__15_\i] {};
  }

  \coordinate (v0__13) at (1, 0, 1);
  \coordinate (v1__13) at (1, 0, 2);
  \coordinate (v2__13) at (1, 1, 1);
  \coordinate (v3__13) at (2, 0, 1);

  \definecolor{vertexcolor__13}{rgb}{ 0 0 0 }

  \tikzstyle{vertexstyle__13_0} = [circle, scale=0.38, fill=vertexcolor__13,]
  \tikzstyle{vertexstyle__13_1} = [circle, scale=0.38, fill=vertexcolor__13,]
  \tikzstyle{vertexstyle__13_2} = [circle, scale=0.38, fill=vertexcolor__13,]
  \tikzstyle{vertexstyle__13_3} = [circle, scale=0.38, fill=vertexcolor__13,]

  \definecolor{facetcolor__13}{rgb}{ 1 1 1 }

  \definecolor{edgecolor__13}{rgb}{ 0 0 0 }
  \tikzstyle{facetstyle__13} = [fill=facetcolor__13, fill opacity=0.25, draw=edgecolor__13, line width=1 pt, line cap=round, line join=round]

  \draw[facetstyle__13] (v0__13) -- (v1__13) -- (v2__13) -- (v0__13) -- cycle;
  \draw[facetstyle__13] (v3__13) -- (v1__13) -- (v0__13) -- (v3__13) -- cycle;
  \draw[facetstyle__13] (v0__13) -- (v2__13) -- (v3__13) -- (v0__13) -- cycle;

   \node at (v0__13) [vertexstyle__13_0] {};

  \draw[facetstyle__13] (v2__13) -- (v1__13) -- (v3__13) -- (v2__13) -- cycle;

  \foreach \i in {1,2,3} {
    \node at (v\i__13) [vertexstyle__13_\i] {};
  }

  \coordinate (v0__9) at (0, 1, 1);
  \coordinate (v1__9) at (0, 1, 2);
  \coordinate (v2__9) at (0, 2, 0);
  \coordinate (v3__9) at (0, 2, 1);
  \coordinate (v4__9) at (1, 1, 1);
  \coordinate (v5__9) at (1, 2, 0);

  \definecolor{vertexcolor__9}{rgb}{ 0 0 0 }

  \tikzstyle{vertexstyle__9_0} = [circle, scale=0.38, fill=vertexcolor__9,]
  \tikzstyle{vertexstyle__9_1} = [circle, scale=0.38, fill=vertexcolor__9,]
  \tikzstyle{vertexstyle__9_2} = [circle, scale=0.38, fill=vertexcolor__9,]
  \tikzstyle{vertexstyle__9_3} = [circle, scale=0.38, fill=vertexcolor__9,]
  \tikzstyle{vertexstyle__9_4} = [circle, scale=0.38, fill=vertexcolor__9,]
  \tikzstyle{vertexstyle__9_5} = [circle, scale=0.38, fill=vertexcolor__9,]

  \definecolor{facetcolor__9}{rgb}{ 0.67843137254 0.84705882352 0.90196078431 }

  \definecolor{edgecolor__9}{rgb}{ 0 0 0 }
  \tikzstyle{facetstyle__9} = [fill=facetcolor__9, fill opacity=0.6, draw=edgecolor__9, line width=1 pt, line cap=round, line join=round]

  \draw[facetstyle__9] (v5__9) -- (v4__9) -- (v0__9) -- (v2__9) -- (v5__9) -- cycle;
  \draw[facetstyle__9] (v2__9) -- (v0__9) -- (v1__9) -- (v3__9) -- (v2__9) -- cycle;
  \draw[facetstyle__9] (v0__9) -- (v4__9) -- (v1__9) -- (v0__9) -- cycle;

   \node at (v0__9) [vertexstyle__9_0] {};

  \draw[facetstyle__9] (v1__9) -- (v4__9) -- (v5__9) -- (v3__9) -- (v1__9) -- cycle;
  \draw[facetstyle__9] (v5__9) -- (v2__9) -- (v3__9) -- (v5__9) -- cycle;

  \foreach \i in {1,3,4,2,5} {
    \node at (v\i__9) [vertexstyle__9_\i] {};
  }

  \coordinate (v0__11) at (0, 2, 0);
  \coordinate (v1__11) at (0, 2, 1);
  \coordinate (v2__11) at (0, 3, 0);
  \coordinate (v3__11) at (1, 2, 0);

  \definecolor{vertexcolor__11}{rgb}{ 0 0 0 }

  \tikzstyle{vertexstyle__11_0} = [circle, scale=0.38, fill=vertexcolor__11,]
  \tikzstyle{vertexstyle__11_1} = [circle, scale=0.38, fill=vertexcolor__11,]
  \tikzstyle{vertexstyle__11_2} = [circle, scale=0.38, fill=vertexcolor__11,]
  \tikzstyle{vertexstyle__11_3} = [circle, scale=0.38, fill=vertexcolor__11,]

  \definecolor{facetcolor__11}{rgb}{ 1 1 1 }

  \definecolor{edgecolor__11}{rgb}{ 0 0 0 }
  \tikzstyle{facetstyle__11} = [fill=facetcolor__11, fill opacity=0.25, draw=edgecolor__11, line width=1 pt, line cap=round, line join=round]

  \draw[facetstyle__11] (v0__11) -- (v1__11) -- (v2__11) -- (v0__11) -- cycle;
  \draw[facetstyle__11] (v3__11) -- (v1__11) -- (v0__11) -- (v3__11) -- cycle;
  \draw[facetstyle__11] (v0__11) -- (v2__11) -- (v3__11) -- (v0__11) -- cycle;

   \node at (v0__11) [vertexstyle__11_0] {};

  \draw[facetstyle__11] (v2__11) -- (v1__11) -- (v3__11) -- (v2__11) -- cycle;

  \foreach \i in {1,2,3} {
    \node at (v\i__11) [vertexstyle__11_\i] {};
  }

  \coordinate (v0__3) at (0, 0, 2);
  \coordinate (v1__3) at (0, 0, 3);
  \coordinate (v2__3) at (0, 1, 2);
  \coordinate (v3__3) at (1, 1, 1);

  \definecolor{vertexcolor__3}{rgb}{ 0 0 0 }

  \tikzstyle{vertexstyle__3_0} = [circle, scale=0.38, fill=vertexcolor__3,]
  \tikzstyle{vertexstyle__3_1} = [circle, scale=0.38, fill=vertexcolor__3,]
  \tikzstyle{vertexstyle__3_2} = [circle, scale=0.38, fill=vertexcolor__3,]
  \tikzstyle{vertexstyle__3_3} = [circle, scale=0.38, fill=vertexcolor__3,]

  \definecolor{facetcolor__3}{rgb}{ 1 1 1 }

  \definecolor{edgecolor__3}{rgb}{ 0 0 0 }
  \tikzstyle{facetstyle__3} = [fill=facetcolor__3, fill opacity=0.25, draw=edgecolor__3, line width=1 pt, line cap=round, line join=round]

  \draw[facetstyle__3] (v0__3) -- (v1__3) -- (v2__3) -- (v0__3) -- cycle;
  \draw[facetstyle__3] (v0__3) -- (v2__3) -- (v3__3) -- (v0__3) -- cycle;
  \draw[facetstyle__3] (v3__3) -- (v1__3) -- (v0__3) -- (v3__3) -- cycle;
  \draw[facetstyle__3] (v2__3) -- (v1__3) -- (v3__3) -- (v2__3) -- cycle;

  \foreach \i in {1,0,2,3} {
    \node at (v\i__3) [vertexstyle__3_\i] {};
  }

  \coordinate (v0__4) at (0, 0, 2);
  \coordinate (v1__4) at (0, 0, 3);
  \coordinate (v2__4) at (1, 0, 2);
  \coordinate (v3__4) at (1, 1, 1);

  \definecolor{vertexcolor__4}{rgb}{ 0 0 0 }

  \tikzstyle{vertexstyle__4_0} = [circle, scale=0.38, fill=vertexcolor__4,]
  \tikzstyle{vertexstyle__4_1} = [circle, scale=0.38, fill=vertexcolor__4,]
  \tikzstyle{vertexstyle__4_2} = [circle, scale=0.38, fill=vertexcolor__4,]
  \tikzstyle{vertexstyle__4_3} = [circle, scale=0.38, fill=vertexcolor__4,]

  \definecolor{facetcolor__4}{rgb}{ 1 1 1 }

  \definecolor{edgecolor__4}{rgb}{ 0 0 0 }
  \tikzstyle{facetstyle__4} = [fill=facetcolor__4, fill opacity=0.25, draw=edgecolor__4, line width=1 pt, line cap=round, line join=round]

  \draw[facetstyle__4] (v1__4) -- (v0__4) -- (v2__4) -- (v1__4) -- cycle;
  \draw[facetstyle__4] (v3__4) -- (v0__4) -- (v1__4) -- (v3__4) -- cycle;
  \draw[facetstyle__4] (v2__4) -- (v0__4) -- (v3__4) -- (v2__4) -- cycle;

   \node at (v0__4) [vertexstyle__4_0] {};

  \draw[facetstyle__4] (v1__4) -- (v2__4) -- (v3__4) -- (v1__4) -- cycle;

  \foreach \i in {1,2,3} {
    \node at (v\i__4) [vertexstyle__4_\i] {};
  }

\end{tikzpicture}\hfill
\begin{tikzpicture}[x  = {(0.9cm,-0.076cm)},
                    y  = {(-0.06cm,0.95cm)},
                    z  = {(-0.44cm,-0.29cm)},
                    scale = 1,
                    color = {lightgray}]

  \node at (0, 0, 3) [] {};

  \coordinate (v0__2) at (0, 0, 1);
  \coordinate (v1__2) at (0, 0, 2);
  \coordinate (v2__2) at (0, 1, 0);
  \coordinate (v3__2) at (0, 1, 1);
  \coordinate (v4__2) at (1, 0, 0);
  \coordinate (v5__2) at (1, 0, 1);

  \definecolor{vertexcolor__2}{rgb}{ 0 0 0 }

  \tikzstyle{vertexstyle__2_0} = [circle, scale=0.38, fill=vertexcolor__2,]
  \tikzstyle{vertexstyle__2_1} = [circle, scale=0.38, fill=vertexcolor__2,]
  \tikzstyle{vertexstyle__2_2} = [circle, scale=0.38, fill=vertexcolor__2,]
  \tikzstyle{vertexstyle__2_3} = [circle, scale=0.38, fill=vertexcolor__2,]
  \tikzstyle{vertexstyle__2_4} = [circle, scale=0.38, fill=vertexcolor__2,]
  \tikzstyle{vertexstyle__2_5} = [circle, scale=0.38, fill=vertexcolor__2,]

  \definecolor{facetcolor__2}{rgb}{ 0.67843137254 0.84705882352 0.90196078431 }

  \definecolor{edgecolor__2}{rgb}{ 0 0 0 }
  \tikzstyle{facetstyle__2} = [fill=facetcolor__2, fill opacity=0.4, draw=edgecolor__2, line width=1 pt, line cap=round, line join=round]

  \draw[facetstyle__2] (v4__2) -- (v5__2) -- (v1__2) -- (v0__2) -- (v4__2) -- cycle;
  \draw[facetstyle__2] (v4__2) -- (v0__2) -- (v2__2) -- (v4__2) -- cycle;
  \draw[facetstyle__2] (v0__2) -- (v1__2) -- (v3__2) -- (v2__2) -- (v0__2) -- cycle;

   \node at (v0__2) [vertexstyle__2_0] {};

  \draw[facetstyle__2] (v3__2) -- (v5__2) -- (v4__2) -- (v2__2) -- (v3__2) -- cycle;
  \draw[facetstyle__2] (v1__2) -- (v5__2) -- (v3__2) -- (v1__2) -- cycle;

  \foreach \i in {1,3,5,2,4} {
    \node at (v\i__2) [vertexstyle__2_\i] {};
  }

  \coordinate (v0__8) at (0, 1, 0);
  \coordinate (v1__8) at (0, 1, 1);
  \coordinate (v2__8) at (0, 2, 0);
  \coordinate (v3__8) at (1, 0, 0);
  \coordinate (v4__8) at (1, 0, 1);
  \coordinate (v5__8) at (1, 1, 0);

  \definecolor{vertexcolor__8}{rgb}{ 0 0 0 }

  \tikzstyle{vertexstyle__8_0} = [circle, scale=0.38, fill=vertexcolor__8,]
  \tikzstyle{vertexstyle__8_1} = [circle, scale=0.38, fill=vertexcolor__8,]
  \tikzstyle{vertexstyle__8_2} = [circle, scale=0.38, fill=vertexcolor__8,]
  \tikzstyle{vertexstyle__8_3} = [circle, scale=0.38, fill=vertexcolor__8,]
  \tikzstyle{vertexstyle__8_4} = [circle, scale=0.38, fill=vertexcolor__8,]
  \tikzstyle{vertexstyle__8_5} = [circle, scale=0.38, fill=vertexcolor__8,]

  \definecolor{facetcolor__8}{rgb}{ 0.67843137254 0.84705882352 0.90196078431 }

  \definecolor{edgecolor__8}{rgb}{ 0 0 0 }
  \tikzstyle{facetstyle__8} = [fill=facetcolor__8, fill opacity=0.4, draw=edgecolor__8, line width=1 pt, line cap=round, line join=round]

  \draw[facetstyle__8] (v2__8) -- (v5__8) -- (v3__8) -- (v0__8) -- (v2__8) -- cycle;
  \draw[facetstyle__8] (v0__8) -- (v3__8) -- (v4__8) -- (v1__8) -- (v0__8) -- cycle;
  \draw[facetstyle__8] (v2__8) -- (v0__8) -- (v1__8) -- (v2__8) -- cycle;

   \node at (v0__8) [vertexstyle__8_0] {};

  \draw[facetstyle__8] (v4__8) -- (v5__8) -- (v2__8) -- (v1__8) -- (v4__8) -- cycle;
  \draw[facetstyle__8] (v3__8) -- (v5__8) -- (v4__8) -- (v3__8) -- cycle;

  \foreach \i in {1,4,3,2,5} {
    \node at (v\i__8) [vertexstyle__8_\i] {};
  }

  \coordinate (v0__10) at (0, 1, 1);
  \coordinate (v1__10) at (0, 2, 0);
  \coordinate (v2__10) at (1, 0, 1);
  \coordinate (v3__10) at (1, 1, 0);
  \coordinate (v4__10) at (1, 1, 1);
  \coordinate (v5__10) at (1, 2, 0);

  \definecolor{vertexcolor__10}{rgb}{ 0 0 0 }

  \tikzstyle{vertexstyle__10_0} = [circle, scale=0.38, fill=vertexcolor__10,]
  \tikzstyle{vertexstyle__10_1} = [circle, scale=0.38, fill=vertexcolor__10,]
  \tikzstyle{vertexstyle__10_2} = [circle, scale=0.38, fill=vertexcolor__10,]
  \tikzstyle{vertexstyle__10_3} = [circle, scale=0.38, fill=vertexcolor__10,]
  \tikzstyle{vertexstyle__10_4} = [circle, scale=0.38, fill=vertexcolor__10,]
  \tikzstyle{vertexstyle__10_5} = [circle, scale=0.38, fill=vertexcolor__10,]

  \definecolor{facetcolor__10}{rgb}{ 0.67843137254 0.84705882352 0.90196078431 }

  \definecolor{edgecolor__10}{rgb}{ 0 0 0 }
  \tikzstyle{facetstyle__10} = [fill=facetcolor__10, fill opacity=0.4, draw=edgecolor__10, line width=1 pt, line cap=round, line join=round]

  \draw[facetstyle__10] (v0__10) -- (v1__10) -- (v3__10) -- (v2__10) -- (v0__10) -- cycle;
  \draw[facetstyle__10] (v1__10) -- (v5__10) -- (v3__10) -- (v1__10) -- cycle;
  \draw[facetstyle__10] (v4__10) -- (v5__10) -- (v1__10) -- (v0__10) -- (v4__10) -- cycle;
  \draw[facetstyle__10] (v3__10) -- (v5__10) -- (v4__10) -- (v2__10) -- (v3__10) -- cycle;
  \draw[facetstyle__10] (v4__10) -- (v0__10) -- (v2__10) -- (v4__10) -- cycle;

  \foreach \i in {0,2,4,1,3,5} {
    \node at (v\i__10) [vertexstyle__10_\i] {};
  }

  \coordinate (v0__14) at (1, 0, 1);
  \coordinate (v1__14) at (1, 1, 0);
  \coordinate (v2__14) at (1, 1, 1);
  \coordinate (v3__14) at (1, 2, 0);
  \coordinate (v4__14) at (2, 0, 0);
  \coordinate (v5__14) at (2, 1, 0);

  \definecolor{vertexcolor__14}{rgb}{ 0 0 0 }

  \tikzstyle{vertexstyle__14_0} = [circle, scale=0.38, fill=vertexcolor__14,]
  \tikzstyle{vertexstyle__14_1} = [circle, scale=0.38, fill=vertexcolor__14,]
  \tikzstyle{vertexstyle__14_2} = [circle, scale=0.38, fill=vertexcolor__14,]
  \tikzstyle{vertexstyle__14_3} = [circle, scale=0.38, fill=vertexcolor__14,]
  \tikzstyle{vertexstyle__14_4} = [circle, scale=0.38, fill=vertexcolor__14,]
  \tikzstyle{vertexstyle__14_5} = [circle, scale=0.38, fill=vertexcolor__14,]

  \definecolor{facetcolor__14}{rgb}{ 0.67843137254 0.84705882352 0.90196078431 }

  \definecolor{edgecolor__14}{rgb}{ 0 0 0 }
  \tikzstyle{facetstyle__14} = [fill=facetcolor__14, fill opacity=0.4, draw=edgecolor__14, line width=1 pt, line cap=round, line join=round]

  \draw[facetstyle__14] (v3__14) -- (v5__14) -- (v4__14) -- (v1__14) -- (v3__14) -- cycle;
  \draw[facetstyle__14] (v0__14) -- (v2__14) -- (v3__14) -- (v1__14) -- (v0__14) -- cycle;
  \draw[facetstyle__14] (v4__14) -- (v0__14) -- (v1__14) -- (v4__14) -- cycle;

   \node at (v1__14) [vertexstyle__14_1] {};

  \draw[facetstyle__14] (v4__14) -- (v5__14) -- (v2__14) -- (v0__14) -- (v4__14) -- cycle;
  \draw[facetstyle__14] (v2__14) -- (v5__14) -- (v3__14) -- (v2__14) -- cycle;

  \foreach \i in {0,2,4,3,5} {
    \node at (v\i__14) [vertexstyle__14_\i] {};
  }

  \coordinate (v0__15) at (1, 0, 1);
  \coordinate (v1__15) at (1, 1, 1);
  \coordinate (v2__15) at (2, 0, 0);
  \coordinate (v3__15) at (2, 0, 1);
  \coordinate (v4__15) at (2, 1, 0);
  \coordinate (v5__15) at (3, 0, 0);

  \definecolor{vertexcolor__15}{rgb}{ 0 0 0 }

  \tikzstyle{vertexstyle__15_0} = [circle, scale=0.38, fill=vertexcolor__15,]
  \tikzstyle{vertexstyle__15_1} = [circle, scale=0.38, fill=vertexcolor__15,]
  \tikzstyle{vertexstyle__15_2} = [circle, scale=0.38, fill=vertexcolor__15,]
  \tikzstyle{vertexstyle__15_3} = [circle, scale=0.38, fill=vertexcolor__15,]
  \tikzstyle{vertexstyle__15_4} = [circle, scale=0.38, fill=vertexcolor__15,]
  \tikzstyle{vertexstyle__15_5} = [circle, scale=0.38, fill=vertexcolor__15,]

  \definecolor{facetcolor__15}{rgb}{ 0.67843137254 0.84705882352 0.90196078431 }

  \definecolor{edgecolor__15}{rgb}{ 0 0 0 }
  \tikzstyle{facetstyle__15} = [fill=facetcolor__15, fill opacity=0.4, draw=edgecolor__15, line width=1 pt, line cap=round, line join=round]

  \draw[facetstyle__15] (v5__15) -- (v3__15) -- (v0__15) -- (v2__15) -- (v5__15) -- cycle;
  \draw[facetstyle__15] (v2__15) -- (v0__15) -- (v1__15) -- (v4__15) -- (v2__15) -- cycle;
  \draw[facetstyle__15] (v5__15) -- (v2__15) -- (v4__15) -- (v5__15) -- cycle;

   \node at (v2__15) [vertexstyle__15_2] {};

  \draw[facetstyle__15] (v1__15) -- (v3__15) -- (v5__15) -- (v4__15) -- (v1__15) -- cycle;
  \draw[facetstyle__15] (v0__15) -- (v3__15) -- (v1__15) -- (v0__15) -- cycle;

  \foreach \i in {0,1,3,4,5} {
    \node at (v\i__15) [vertexstyle__15_\i] {};
  }

  \coordinate (v0__9) at (0, 1, 1);
  \coordinate (v1__9) at (0, 1, 2);
  \coordinate (v2__9) at (0, 2, 0);
  \coordinate (v3__9) at (0, 2, 1);
  \coordinate (v4__9) at (1, 1, 1);
  \coordinate (v5__9) at (1, 2, 0);

  \definecolor{vertexcolor__9}{rgb}{ 0 0 0 }

  \tikzstyle{vertexstyle__9_0} = [circle, scale=0.38, fill=vertexcolor__9,]
  \tikzstyle{vertexstyle__9_1} = [circle, scale=0.38, fill=vertexcolor__9,]
  \tikzstyle{vertexstyle__9_2} = [circle, scale=0.38, fill=vertexcolor__9,]
  \tikzstyle{vertexstyle__9_3} = [circle, scale=0.38, fill=vertexcolor__9,]
  \tikzstyle{vertexstyle__9_4} = [circle, scale=0.38, fill=vertexcolor__9,]
  \tikzstyle{vertexstyle__9_5} = [circle, scale=0.38, fill=vertexcolor__9,]

  \definecolor{facetcolor__9}{rgb}{ 0.67843137254 0.84705882352 0.90196078431 }

  \definecolor{edgecolor__9}{rgb}{ 0 0 0 }
  \tikzstyle{facetstyle__9} = [fill=facetcolor__9, fill opacity=0.4, draw=edgecolor__9, line width=1 pt, line cap=round, line join=round]

  \draw[facetstyle__9] (v5__9) -- (v4__9) -- (v0__9) -- (v2__9) -- (v5__9) -- cycle;
  \draw[facetstyle__9] (v2__9) -- (v0__9) -- (v1__9) -- (v3__9) -- (v2__9) -- cycle;
  \draw[facetstyle__9] (v0__9) -- (v4__9) -- (v1__9) -- (v0__9) -- cycle;

   \node at (v0__9) [vertexstyle__9_0] {};

  \draw[facetstyle__9] (v1__9) -- (v4__9) -- (v5__9) -- (v3__9) -- (v1__9) -- cycle;
  \draw[facetstyle__9] (v5__9) -- (v2__9) -- (v3__9) -- (v5__9) -- cycle;

  \foreach \i in {1,3,4,2,5} {
    \node at (v\i__9) [vertexstyle__9_\i] {};
  }

\end{tikzpicture}\hfill
\begin{tikzpicture}[scale=1]
\draw (0,1) -- (0,0);
\draw (0,0) -- (1.73205080756888,-1);
\draw (0,0) -- (-1.73205080756888,-1);
\fill (1.73205080756888, -1) circle (2.85pt);
\fill (0.866025403784439, -0.5) circle (2.85pt);
\fill (0, 0) circle (2.85pt);
\fill (-0.866025403784438, -0.5) circle (2.85pt);
\fill (-1.73205080756888, -1) circle (2.85pt);
\fill (0, 1) circle (2.85pt);
\node at (-2.3, -1.7) {};
\end{tikzpicture}

%% file: 3.tex
\begin{tikzpicture}[scale=0.18]
\draw (3.48625962311572,-1.26889473173823) -- (2.9,0);
\draw (4.439,0) -- (3.48625962311572,1.26889473173823);
\draw (-0.644234739144311,3.65363676367529) -- (-1.45,2.51147367097487);
\draw (2.9,0) -- (2,0);
\draw (2.9,0) -- (3.48625962311572,1.26889473173823);
\draw (-1.45,2.51147367097487) -- (-1,1.73205080756888);
\draw (-1.45,2.51147367097487) -- (-2.84202488397141,2.38474203193706);
\draw (2,0) -- (1,0);
\draw (-1,1.73205080756888) -- (-0.5,0.866025403784439);
\draw (1,0) -- (-0.5,0.866025403784439);
\draw (1,0) -- (-0.5,-0.866025403784438);
\draw (-0.5,0.866025403784439) -- (-0.5,-0.866025403784438);
\draw (-2.2195,3.84428676739912) -- (-2.84202488397141,2.38474203193706);
\draw (-0.5,-0.866025403784438) -- (-1,-1.73205080756888);
\draw (-1,-1.73205080756888) -- (-1.45,-2.51147367097487);
\draw (-1.855,-3.21295424804027) -- (-1.45,-2.51147367097487);
\fill[black] (-2.2195, 3.84428676739912) circle (5.43130215351701pt);
\fill[black] (-0.644234739144311, 3.65363676367529) circle (5.43130215351701pt);
\fill[black] (4.439, 0) circle (5.43130215351701pt);
\fill[black] (-1.45, 2.51147367097487) circle (5.43130215351701pt);
\fill[black] (-0.5, 0.866025403784439) circle (5.43130215351701pt);
\fill[black] (-1, 1.73205080756888) circle (5.43130215351701pt);
\fill[black] (-0.5, -0.866025403784438) circle (5.43130215351701pt);
\fill[black] (3.48625962311572, 1.26889473173823) circle (5.43130215351701pt);
\fill[black] (3.48625962311572, -1.26889473173823) circle (5.43130215351701pt);
\fill[black] (1, 0) circle (5.43130215351701pt);
\fill[black] (-1.855, -3.21295424804027) circle (5.43130215351701pt);
\fill[black] (-2.84202488397141, 2.38474203193706) circle (5.43130215351701pt);
\fill[black] (-1, -1.73205080756888) circle (5.43130215351701pt);
\fill[black] (2, 0) circle (5.43130215351701pt);
\fill[black] (-1.45, -2.51147367097487) circle (5.43130215351701pt);
\fill[black] (2.9, 0) circle (5.43130215351701pt);
\end{tikzpicture}

%% file: 4.tex
\begin{tikzpicture}[scale=0.21]
\draw (1.22464679914735e-16,2) -- (6.12323399573677e-17,1);
\draw (1.22464679914735e-16,2) -- (1.77573785876366e-16,2.9);
\draw (2.8011848962383,-0.75057523079731) -- (2,0);
\draw (2.8011848962383,-0.75057523079731) -- (3.58358481553244,-0.960218657330352);
\draw (2,0) -- (1,0);
\draw (2,0) -- (2.8011848962383,0.75057523079731);
\draw (1,0) -- (6.12323399573677e-17,1);
\draw (1,0) -- (-1.83697019872103e-16,-1);
\draw (6.12323399573677e-17,1) -- (-1,1.22464679914735e-16);
\draw (0.960218657330353,3.58358481553244) -- (1.77573785876366e-16,2.9);
\draw (1.77573785876366e-16,2.9) -- (-0.960218657330351,3.58358481553244);
\draw (2.8011848962383,0.75057523079731) -- (3.58358481553244,0.960218657330352);
\draw (-1,1.22464679914735e-16) -- (-1.83697019872103e-16,-1);
\draw (3.58358481553244,0.960218657330352) -- (4.439,0);
\draw (-0.960218657330351,3.58358481553244) -- (1.14889774121009,4.28774474289717);
\draw (-0.960218657330351,3.58358481553244) -- (-1.14889774121009,4.28774474289717);
\fill[black] (-1.14889774121009, 4.28774474289717) circle (4.75pt);
\fill[black] (1, 0) circle (4.75pt);
\fill[black] (1.14889774121009, 4.28774474289717) circle (4.75pt);
\fill[black] (3.58358481553244, -0.960218657330352) circle (4.75pt);
\fill[black] (-0.960218657330351, 3.58358481553244) circle (4.75pt);
\fill[black] (0.960218657330353, 3.58358481553244) circle (4.75pt);
\fill[black] (4.439, 0) circle (4.75pt);
\fill[black] (1.77573785876366e-16, 2.9) circle (4.75pt);
\fill[black] (1.22464679914735e-16, 2) circle (4.75pt);
\fill[black] (-1, 1.22464679914735e-16) circle (4.75pt);
\fill[black] (-1.83697019872103e-16, -1) circle (4.75pt);
\fill[black] (2.8011848962383, -0.75057523079731) circle (4.75pt);
\fill[black] (3.58358481553244, 0.960218657330352) circle (4.75pt);
\fill[black] (2, 0) circle (4.75pt);
\fill[black] (2.8011848962383, 0.75057523079731) circle (4.75pt);
\fill[black] (6.12323399573677e-17, 1) circle (4.75pt);
\end{tikzpicture}

%% file: 5.tex
\begin{tikzpicture}[scale=0.22]
\draw (2.75806389725595,0.896149283687348) -- (1.90211303259031,0.618033988749895);
\draw (2.75806389725595,0.896149283687348) -- (3.52841967545502,1.14645304913105);
\draw (1.90211303259031,0.618033988749895) -- (0.951056516295154,0.309016994374947);
\draw (0,2.9) -- (0,2);
\draw (0,2.9) -- (0,3.71);
\draw (0,2) -- (0,1);
\draw (0,3.71) -- (0,4.439);
\draw (3.52841967545502,1.14645304913105) -- (4.22173987583419,1.37172643803039);
\draw (0,1) -- (0.951056516295154,0.309016994374947);
\draw (0,1) -- (-0.951056516295154,0.309016994374947);
\draw (0.951056516295154,0.309016994374947) -- (0.587785252292473,-0.809016994374947);
\draw (4.84572805617544,1.57447248803979) -- (4.22173987583419,1.37172643803039);
\draw (0.587785252292473,-0.809016994374947) -- (-0.587785252292473,-0.809016994374948);
\draw (-0.951056516295154,0.309016994374947) -- (-0.587785252292473,-0.809016994374948);
\draw (-0.951056516295154,0.309016994374947) -- (-1.90211303259031,0.618033988749894);
\draw (-1.90211303259031,0.618033988749894) -- (-2.75806389725595,0.896149283687347);
\fill[black] (-0.587785252292473, -0.809016994374948) circle (4.26520876399966pt);
\fill[black] (0, 4.439) circle (4.26520876399966pt);
\fill[black] (0.587785252292473, -0.809016994374947) circle (4.26520876399966pt);
\fill[black] (0, 1) circle (4.26520876399966pt);
\fill[black] (-1.90211303259031, 0.618033988749894) circle (4.26520876399966pt);
\fill[black] (0, 2) circle (4.26520876399966pt);
\fill[black] (-2.75806389725595, 0.896149283687347) circle (4.26520876399966pt);
\fill[black] (0, 3.71) circle (4.26520876399966pt);
\fill[black] (0.951056516295154, 0.309016994374947) circle (4.26520876399966pt);
\fill[black] (1.90211303259031, 0.618033988749895) circle (4.26520876399966pt);
\fill[black] (0, 2.9) circle (4.26520876399966pt);
\fill[black] (4.22173987583419, 1.37172643803039) circle (4.26520876399966pt);
\fill[black] (4.84572805617544, 1.57447248803979) circle (4.26520876399966pt);
\fill[black] (2.75806389725595, 0.896149283687348) circle (4.26520876399966pt);
\fill[black] (-0.951056516295154, 0.309016994374947) circle (4.26520876399966pt);
\fill[black] (3.52841967545502, 1.14645304913105) circle (4.26520876399966pt);
\end{tikzpicture}

%% file: 6.tex
\begin{tikzpicture}[scale=0.26]
\draw (2.84202488397141,2.38474203193706) -- (2.22152888504504,1.86408406809096);
\draw (0,1) -- (0.866025403784439,0.5);
\draw (0,1) -- (-0.866025403784439,0.5);
\draw (0.866025403784439,0.5) -- (0.866025403784439,-0.5);
\draw (0.866025403784439,0.5) -- (1.73205080756888,1);
\draw (0.866025403784439,-0.5) -- (1.22464679914735e-16,-1);
\draw (-0.866025403784439,0.5) -- (-0.866025403784438,-0.5);
\draw (-0.866025403784439,0.5) -- (-1.73205080756888,1);
\draw (-0.866025403784438,-0.5) -- (1.22464679914735e-16,-1);
\draw (-2.72510860027913,0.991858415644438) -- (-1.73205080756888,1);
\draw (-1.73205080756888,1) -- (-2.22152888504504,1.86408406809096);
\draw (2.22152888504504,1.86408406809096) -- (1.73205080756888,1);
\draw (1.73205080756888,1) -- (2.72510860027913,0.99185841564444);
\draw (3.84428676739912,2.2195) -- (3.48625962311572,1.26889473173823);
\draw (3.48625962311572,1.26889473173823) -- (2.72510860027913,0.99185841564444);
\draw (1.22464679914735e-16,-1) -- (2.44929359829471e-16,-2);
\fill[black] (-2.22152888504504, 1.86408406809096) circle (3.8909863237109pt);
\fill[black] (-0.866025403784439, 0.5) circle (3.8909863237109pt);
\fill[black] (0.866025403784439, 0.5) circle (3.8909863237109pt);
\fill[black] (1.73205080756888, 1) circle (3.8909863237109pt);
\fill[black] (0, 1) circle (3.8909863237109pt);
\fill[black] (3.48625962311572, 1.26889473173823) circle (3.8909863237109pt);
\fill[black] (2.84202488397141, 2.38474203193706) circle (3.8909863237109pt);
\fill[black] (2.22152888504504, 1.86408406809096) circle (3.8909863237109pt);
\fill[black] (-2.72510860027913, 0.991858415644438) circle (3.8909863237109pt);
\fill[black] (3.84428676739912, 2.2195) circle (3.8909863237109pt);
\fill[black] (2.72510860027913, 0.99185841564444) circle (3.8909863237109pt);
\fill[black] (-0.866025403784438, -0.5) circle (3.8909863237109pt);
\fill[black] (1.22464679914735e-16, -1) circle (3.8909863237109pt);
\fill[black] (-1.73205080756888, 1) circle (3.8909863237109pt);
\fill[black] (0.866025403784439, -0.5) circle (3.8909863237109pt);
\fill[black] (2.44929359829471e-16, -2) circle (3.8909863237109pt);
\end{tikzpicture}

%% file: 7.tex
\begin{tikzpicture}[scale=0.27]
\draw (0,1) -- (0.78183148246803,0.623489801858734);
\draw (0,1) -- (-0.78183148246803,0.623489801858733);
\draw (0.78183148246803,0.623489801858734) -- (0.974927912181824,-0.222520933956314);
\draw (-0.867767478235116,-1.80193773580484) -- (-0.433883739117558,-0.900968867902419);
\draw (-0.78183148246803,0.623489801858733) -- (-1.56366296493606,1.24697960371747);
\draw (-0.78183148246803,0.623489801858733) -- (-0.974927912181824,-0.222520933956315);
\draw (-1.56366296493606,1.24697960371747) -- (-2.26731129915729,1.80812042539033);
\draw (-0.974927912181824,-0.222520933956315) -- (-0.433883739117558,-0.900968867902419);
\draw (-0.974927912181824,-0.222520933956315) -- (-1.94985582436365,-0.445041867912629);
\draw (-0.433883739117558,-0.900968867902419) -- (0.433883739117558,-0.900968867902419);
\draw (-2.26731129915729,1.80812042539033) -- (-2.90059479995639,2.3131471648959);
\draw (2.89189101182542,-0.21671727140063) -- (1.94985582436365,-0.445041867912629);
\draw (0.974927912181824,-0.222520933956314) -- (1.94985582436365,-0.445041867912629);
\draw (0.974927912181824,-0.222520933956314) -- (0.433883739117558,-0.900968867902419);
\draw (1.94985582436365,-0.445041867912629) -- (2.69953387106819,-1.05948897066254);
\draw (2.69953387106819,-1.05948897066254) -- (3.61698255419457,-0.825552664977926);
\fill[black] (-1.56366296493606, 1.24697960371747) circle (3.58621578407382pt);
\fill[black] (-2.26731129915729, 1.80812042539033) circle (3.58621578407382pt);
\fill[black] (2.89189101182542, -0.21671727140063) circle (3.58621578407382pt);
\fill[black] (-0.867767478235116, -1.80193773580484) circle (3.58621578407382pt);
\fill[black] (0.78183148246803, 0.623489801858734) circle (3.58621578407382pt);
\fill[black] (1.94985582436365, -0.445041867912629) circle (3.58621578407382pt);
\fill[black] (-2.90059479995639, 2.3131471648959) circle (3.58621578407382pt);
\fill[black] (0, 1) circle (3.58621578407382pt);
\fill[black] (-0.433883739117558, -0.900968867902419) circle (3.58621578407382pt);
\fill[black] (0.974927912181824, -0.222520933956314) circle (3.58621578407382pt);
\fill[black] (0.433883739117558, -0.900968867902419) circle (3.58621578407382pt);
\fill[black] (-0.974927912181824, -0.222520933956315) circle (3.58621578407382pt);
\fill[black] (-1.94985582436365, -0.445041867912629) circle (3.58621578407382pt);
\fill[black] (2.69953387106819, -1.05948897066254) circle (3.58621578407382pt);
\fill[black] (-0.78183148246803, 0.623489801858733) circle (3.58621578407382pt);
\fill[black] (3.61698255419457, -0.825552664977926) circle (3.58621578407382pt);
\end{tikzpicture}

%% file: 8.tex
\begin{tikzpicture}[scale=0.26]
\draw (2.9,0) -- (3.71,0);
\draw (2.9,0) -- (2,0);
\draw (3.71,0) -- (4.439,0);
\draw (2,0) -- (1,0);
\draw (1,0) -- (0.707106781186548,0.707106781186547);
\draw (1,0) -- (0.707106781186547,-0.707106781186548);
\draw (0.707106781186548,0.707106781186547) -- (6.12323399573677e-17,1);
\draw (0.707106781186547,-0.707106781186548) -- (-1.83697019872103e-16,-1);
\draw (-1.83697019872103e-16,-1) -- (-0.707106781186548,-0.707106781186547);
\draw (4.439,0) -- (5.0951,0);
\draw (5.0951,0) -- (5.68559,0);
\draw (6.12323399573677e-17,1) -- (-0.707106781186547,0.707106781186548);
\draw (-0.707106781186548,-0.707106781186547) -- (-1,1.22464679914735e-16);
\draw (-0.707106781186547,0.707106781186548) -- (-1.41421356237309,1.4142135623731);
\draw (-0.707106781186547,0.707106781186548) -- (-1,1.22464679914735e-16);
\draw (-1.41421356237309,1.4142135623731) -- (-2.05060966544099,2.05060966544099);
\fill[black] (1, 0) circle (3.32842712474619pt);
\fill[black] (-2.05060966544099, 2.05060966544099) circle (3.32842712474619pt);
\fill[black] (0.707106781186547, -0.707106781186548) circle (3.32842712474619pt);
\fill[black] (-1.41421356237309, 1.4142135623731) circle (3.32842712474619pt);
\fill[black] (-0.707106781186548, -0.707106781186547) circle (3.32842712474619pt);
\fill[black] (-1, 1.22464679914735e-16) circle (3.32842712474619pt);
\fill[black] (4.439, 0) circle (3.32842712474619pt);
\fill[black] (2.9, 0) circle (3.32842712474619pt);
\fill[black] (5.68559, 0) circle (3.32842712474619pt);
\fill[black] (-0.707106781186547, 0.707106781186548) circle (3.32842712474619pt);
\fill[black] (-1.83697019872103e-16, -1) circle (3.32842712474619pt);
\fill[black] (5.0951, 0) circle (3.32842712474619pt);
\fill[black] (0.707106781186548, 0.707106781186547) circle (3.32842712474619pt);
\fill[black] (3.71, 0) circle (3.32842712474619pt);
\fill[black] (6.12323399573677e-17, 1) circle (3.32842712474619pt);
\fill[black] (2, 0) circle (3.32842712474619pt);
\end{tikzpicture}

%% file: 9.tex
\begin{tikzpicture}[scale=0.3]
\draw (0,1) -- (0.642787609686539,0.766044443118978);
\draw (0,1) -- (-0.64278760968654,0.766044443118978);
\draw (0.642787609686539,0.766044443118978) -- (0.984807753012208,0.17364817766693);
\draw (-0.64278760968654,0.766044443118978) -- (-0.984807753012208,0.17364817766693);
\draw (-1.96961550602442,0.34729635533386) -- (-0.984807753012208,0.17364817766693);
\draw (-1.96961550602442,0.34729635533386) -- (-2.8559424837354,0.503579715234097);
\draw (-0.984807753012208,0.17364817766693) -- (-0.866025403784438,-0.5);
\draw (0.984807753012208,0.17364817766693) -- (0.866025403784439,-0.5);
\draw (0.866025403784439,-0.5) -- (0.342020143325669,-0.939692620785908);
\draw (0.99185841564444,-2.72510860027913) -- (0.684040286651338,-1.87938524157182);
\draw (0.684040286651338,-1.87938524157182) -- (0.342020143325669,-0.939692620785908);
\draw (0.342020143325669,-0.939692620785908) -- (-0.342020143325669,-0.939692620785908);
\draw (-0.866025403784438,-0.5) -- (-1.73205080756888,-1);
\draw (-0.866025403784438,-0.5) -- (-0.342020143325669,-0.939692620785908);
\draw (-1.73205080756888,-1) -- (-2.51147367097487,-1.45);
\draw (-2.51147367097487,-1.45) -- (-3.21295424804027,-1.855);
\fill[black] (-0.64278760968654, 0.766044443118978) circle (3.10416666666667pt);
\fill[black] (0.342020143325669, -0.939692620785908) circle (3.10416666666667pt);
\fill[black] (0.642787609686539, 0.766044443118978) circle (3.10416666666667pt);
\fill[black] (-0.984807753012208, 0.17364817766693) circle (3.10416666666667pt);
\fill[black] (0.99185841564444, -2.72510860027913) circle (3.10416666666667pt);
\fill[black] (0.984807753012208, 0.17364817766693) circle (3.10416666666667pt);
\fill[black] (-1.73205080756888, -1) circle (3.10416666666667pt);
\fill[black] (0.684040286651338, -1.87938524157182) circle (3.10416666666667pt);
\fill[black] (0, 1) circle (3.10416666666667pt);
\fill[black] (0.866025403784439, -0.5) circle (3.10416666666667pt);
\fill[black] (-2.51147367097487, -1.45) circle (3.10416666666667pt);
\fill[black] (-0.866025403784438, -0.5) circle (3.10416666666667pt);
\fill[black] (-2.8559424837354, 0.503579715234097) circle (3.10416666666667pt);
\fill[black] (-0.342020143325669, -0.939692620785908) circle (3.10416666666667pt);
\fill[black] (-3.21295424804027, -1.855) circle (3.10416666666667pt);
\fill[black] (-1.96961550602442, 0.34729635533386) circle (3.10416666666667pt);
\end{tikzpicture}

%% file: 10.tex
\begin{tikzpicture}[scale=0.3]
\draw (0,1) -- (0.587785252292473,0.809016994374947);
\draw (0,1) -- (0,2);
\draw (0,1) -- (-0.587785252292473,0.809016994374947);
\draw (0.587785252292473,0.809016994374947) -- (0.951056516295154,0.309016994374947);
\draw (-0.587785252292473,0.809016994374947) -- (-0.951056516295154,0.309016994374947);
\draw (-0.587785252292473,-0.809016994374948) -- (1.22464679914735e-16,-1);
\draw (-0.587785252292473,-0.809016994374948) -- (-0.951056516295154,-0.309016994374948);
\draw (2.44929359829471e-16,-2) -- (1.22464679914735e-16,-1);
\draw (1.22464679914735e-16,-1) -- (0.587785252292473,-0.809016994374947);
\draw (-0.951056516295154,0.309016994374947) -- (-0.951056516295154,-0.309016994374948);
\draw (-0.951056516295154,0.309016994374947) -- (-1.90211303259031,0.618033988749894);
\draw (0.951056516295154,0.309016994374947) -- (1.90211303259031,0.618033988749895);
\draw (0.951056516295154,0.309016994374947) -- (0.951056516295154,-0.309016994374947);
\draw (1.90211303259031,0.618033988749895) -- (2.75806389725595,0.896149283687348);
\draw (0.951056516295154,-0.309016994374947) -- (0.587785252292473,-0.809016994374947);
\draw (2.75806389725595,0.896149283687348) -- (3.52841967545502,1.14645304913105);
\fill[black] (3.52841967545502, 1.14645304913105) circle (2.9048221281347pt);
\fill[black] (-0.587785252292473, 0.809016994374947) circle (2.9048221281347pt);
\fill[black] (2.75806389725595, 0.896149283687348) circle (2.9048221281347pt);
\fill[black] (-0.951056516295154, 0.309016994374947) circle (2.9048221281347pt);
\fill[black] (1.90211303259031, 0.618033988749895) circle (2.9048221281347pt);
\fill[black] (0.587785252292473, -0.809016994374947) circle (2.9048221281347pt);
\fill[black] (2.44929359829471e-16, -2) circle (2.9048221281347pt);
\fill[black] (1.22464679914735e-16, -1) circle (2.9048221281347pt);
\fill[black] (0.951056516295154, -0.309016994374947) circle (2.9048221281347pt);
\fill[black] (0, 1) circle (2.9048221281347pt);
\fill[black] (-1.90211303259031, 0.618033988749894) circle (2.9048221281347pt);
\fill[black] (0, 2) circle (2.9048221281347pt);
\fill[black] (0.951056516295154, 0.309016994374947) circle (2.9048221281347pt);
\fill[black] (0.587785252292473, 0.809016994374947) circle (2.9048221281347pt);
\fill[black] (-0.951056516295154, -0.309016994374948) circle (2.9048221281347pt);
\fill[black] (-0.587785252292473, -0.809016994374948) circle (2.9048221281347pt);
\end{tikzpicture}

%% file: 11.tex
\begin{tikzpicture}[scale=0.35]
\draw (-0.563465113682859,-1.91898594722899) -- (-0.281732556841429,-0.959492973614497);
\draw (1.0812816349112,1.68250706566236) -- (0.540640817455598,0.841253532831181);
\draw (0,1) -- (0.540640817455598,0.841253532831181);
\draw (0,1) -- (-0.540640817455597,0.841253532831181);
\draw (0.540640817455598,0.841253532831181) -- (0.909631995354518,0.415415013001886);
\draw (-0.755749574354258,-0.654860733945285) -- (-0.281732556841429,-0.959492973614497);
\draw (-0.755749574354258,-0.654860733945285) -- (-0.989821441880933,-0.142314838273285);
\draw (-0.281732556841429,-0.959492973614497) -- (0.28173255684143,-0.959492973614497);
\draw (-0.540640817455597,0.841253532831181) -- (-0.909631995354519,0.415415013001886);
\draw (-0.989821441880933,-0.142314838273285) -- (-0.909631995354519,0.415415013001886);
\draw (-0.909631995354519,0.415415013001886) -- (-1.81926399070904,0.830830026003772);
\draw (0.28173255684143,-0.959492973614497) -- (0.755749574354258,-0.654860733945285);
\draw (0.755749574354258,-0.654860733945285) -- (1.51149914870852,-1.30972146789057);
\draw (0.755749574354258,-0.654860733945285) -- (0.989821441880933,-0.142314838273285);
\draw (0.909631995354518,0.415415013001886) -- (1.81926399070904,0.830830026003773);
\draw (0.909631995354518,0.415415013001886) -- (0.989821441880933,-0.142314838273285);
\fill[black] (0.540640817455598, 0.841253532831181) circle (2.72459075662211pt);
\fill[black] (0.989821441880933, -0.142314838273285) circle (2.72459075662211pt);
\fill[black] (0.755749574354258, -0.654860733945285) circle (2.72459075662211pt);
\fill[black] (0.909631995354518, 0.415415013001886) circle (2.72459075662211pt);
\fill[black] (-0.281732556841429, -0.959492973614497) circle (2.72459075662211pt);
\fill[black] (1.81926399070904, 0.830830026003773) circle (2.72459075662211pt);
\fill[black] (-0.989821441880933, -0.142314838273285) circle (2.72459075662211pt);
\fill[black] (1.51149914870852, -1.30972146789057) circle (2.72459075662211pt);
\fill[black] (-0.563465113682859, -1.91898594722899) circle (2.72459075662211pt);
\fill[black] (-1.81926399070904, 0.830830026003772) circle (2.72459075662211pt);
\fill[black] (-0.540640817455597, 0.841253532831181) circle (2.72459075662211pt);
\fill[black] (-0.909631995354519, 0.415415013001886) circle (2.72459075662211pt);
\fill[black] (-0.755749574354258, -0.654860733945285) circle (2.72459075662211pt);
\fill[black] (0.28173255684143, -0.959492973614497) circle (2.72459075662211pt);
\fill[black] (0, 1) circle (2.72459075662211pt);
\fill[black] (1.0812816349112, 1.68250706566236) circle (2.72459075662211pt);
\end{tikzpicture}

%% file: 12.tex
\begin{tikzpicture}[scale=0.4]
\draw (2.44929359829471e-16,-2) -- (1.22464679914735e-16,-1);
\draw (0,1) -- (0.5,0.866025403784439);
\draw (0,1) -- (-0.5,0.866025403784438);
\draw (0.5,0.866025403784439) -- (0.866025403784439,0.5);
\draw (0.866025403784439,0.5) -- (1,6.12323399573677e-17);
\draw (-0.5,-0.866025403784439) -- (1.22464679914735e-16,-1);
\draw (-0.5,-0.866025403784439) -- (-0.866025403784438,-0.5);
\draw (1.22464679914735e-16,-1) -- (0.5,-0.866025403784439);
\draw (-0.866025403784438,-0.5) -- (-1,-1.83697019872103e-16);
\draw (-0.5,0.866025403784438) -- (-0.866025403784439,0.5);
\draw (-0.866025403784439,0.5) -- (-1,-1.83697019872103e-16);
\draw (-1,-1.83697019872103e-16) -- (-2,-3.67394039744206e-16);
\draw (0.5,-0.866025403784439) -- (1,-1.73205080756888);
\draw (0.5,-0.866025403784439) -- (0.866025403784439,-0.5);
\draw (0.866025403784439,-0.5) -- (1,6.12323399573677e-17);
\draw (0.866025403784439,-0.5) -- (1.73205080756888,-1);
\fill[black] (1.22464679914735e-16, -1) circle (2.5594010767585pt);
\fill[black] (0.866025403784439, -0.5) circle (2.5594010767585pt);
\fill[black] (0.5, -0.866025403784439) circle (2.5594010767585pt);
\fill[black] (1, 6.12323399573677e-17) circle (2.5594010767585pt);
\fill[black] (-0.5, 0.866025403784438) circle (2.5594010767585pt);
\fill[black] (0.866025403784439, 0.5) circle (2.5594010767585pt);
\fill[black] (1.73205080756888, -1) circle (2.5594010767585pt);
\fill[black] (-0.866025403784439, 0.5) circle (2.5594010767585pt);
\fill[black] (-0.5, -0.866025403784439) circle (2.5594010767585pt);
\fill[black] (0, 1) circle (2.5594010767585pt);
\fill[black] (0.5, 0.866025403784439) circle (2.5594010767585pt);
\fill[black] (-2, -3.67394039744206e-16) circle (2.5594010767585pt);
\fill[black] (2.44929359829471e-16, -2) circle (2.5594010767585pt);
\fill[black] (-1, -1.83697019872103e-16) circle (2.5594010767585pt);
\fill[black] (1, -1.73205080756888) circle (2.5594010767585pt);
\fill[black] (-0.866025403784438, -0.5) circle (2.5594010767585pt);
\end{tikzpicture}

%% file: 13.tex
\begin{tikzpicture}[scale=0.4]
\draw (0,1) -- (0.464723172043769,0.88545602565321);
\draw (0,1) -- (-0.464723172043768,0.88545602565321);
\draw (0.464723172043769,0.88545602565321) -- (0.822983865893656,0.568064746731156);
\draw (-0.464723172043768,0.88545602565321) -- (-0.822983865893657,0.568064746731155);
\draw (-0.478631328575115,-1.9418836348521) -- (-0.239315664287557,-0.970941817426052);
\draw (-0.663122658240795,-0.748510748171101) -- (-0.239315664287557,-0.970941817426052);
\draw (-0.663122658240795,-0.748510748171101) -- (-0.935016242685415,-0.354604887042536);
\draw (-0.239315664287557,-0.970941817426052) -- (0.239315664287558,-0.970941817426052);
\draw (-0.822983865893657,0.568064746731155) -- (-0.992708874098054,0.120536680255323);
\draw (-0.992708874098054,0.120536680255323) -- (-0.935016242685415,-0.354604887042536);
\draw (-0.992708874098054,0.120536680255323) -- (-1.98541774819611,0.241073360510646);
\draw (0.822983865893656,0.568064746731156) -- (0.992708874098054,0.120536680255323);
\draw (0.478631328575115,-1.9418836348521) -- (0.239315664287558,-0.970941817426052);
\draw (0.239315664287558,-0.970941817426052) -- (0.663122658240795,-0.748510748171101);
\draw (0.992708874098054,0.120536680255323) -- (0.935016242685415,-0.354604887042535);
\draw (0.935016242685415,-0.354604887042535) -- (0.663122658240795,-0.748510748171101);
\fill[black] (-0.822983865893657, 0.568064746731155) circle (2.40630078490092pt);
\fill[black] (0, 1) circle (2.40630078490092pt);
\fill[black] (-1.98541774819611, 0.241073360510646) circle (2.40630078490092pt);
\fill[black] (0.239315664287558, -0.970941817426052) circle (2.40630078490092pt);
\fill[black] (0.464723172043769, 0.88545602565321) circle (2.40630078490092pt);
\fill[black] (0.822983865893656, 0.568064746731156) circle (2.40630078490092pt);
\fill[black] (-0.464723172043768, 0.88545602565321) circle (2.40630078490092pt);
\fill[black] (-0.935016242685415, -0.354604887042536) circle (2.40630078490092pt);
\fill[black] (-0.663122658240795, -0.748510748171101) circle (2.40630078490092pt);
\fill[black] (0.663122658240795, -0.748510748171101) circle (2.40630078490092pt);
\fill[black] (-0.478631328575115, -1.9418836348521) circle (2.40630078490092pt);
\fill[black] (0.935016242685415, -0.354604887042535) circle (2.40630078490092pt);
\fill[black] (-0.992708874098054, 0.120536680255323) circle (2.40630078490092pt);
\fill[black] (-0.239315664287557, -0.970941817426052) circle (2.40630078490092pt);
\fill[black] (0.992708874098054, 0.120536680255323) circle (2.40630078490092pt);
\fill[black] (0.478631328575115, -1.9418836348521) circle (2.40630078490092pt);
\end{tikzpicture}

%% file: 14.tex
\begin{tikzpicture}[scale=0.4]
\draw (1,0) -- (0.900968867902419,0.433883739117558);
\draw (1,0) -- (2,0);
\draw (1,0) -- (0.900968867902419,-0.433883739117558);
\draw (0.900968867902419,0.433883739117558) -- (0.623489801858734,0.78183148246803);
\draw (0.900968867902419,-0.433883739117558) -- (0.623489801858733,-0.78183148246803);
\draw (-0.900968867902419,-0.433883739117558) -- (-1,1.22464679914735e-16);
\draw (-0.900968867902419,-0.433883739117558) -- (-0.623489801858734,-0.78183148246803);
\draw (-1,1.22464679914735e-16) -- (-0.900968867902419,0.433883739117558);
\draw (-0.623489801858734,-0.78183148246803) -- (-0.222520933956315,-0.974927912181824);
\draw (0.222520933956313,-0.974927912181824) -- (-0.222520933956315,-0.974927912181824);
\draw (0.222520933956313,-0.974927912181824) -- (0.623489801858733,-0.78183148246803);
\draw (0.623489801858734,0.78183148246803) -- (0.222520933956314,0.974927912181824);
\draw (-0.900968867902419,0.433883739117558) -- (-0.623489801858733,0.78183148246803);
\draw (-0.623489801858733,0.78183148246803) -- (-1.24697960371747,1.56366296493606);
\draw (-0.623489801858733,0.78183148246803) -- (-0.222520933956314,0.974927912181824);
\draw (0.222520933956314,0.974927912181824) -- (-0.222520933956314,0.974927912181824);
\fill[black] (-0.623489801858734, -0.78183148246803) circle (2.2630899352994pt);
\fill[black] (-0.623489801858733, 0.78183148246803) circle (2.2630899352994pt);
\fill[black] (1, 0) circle (2.2630899352994pt);
\fill[black] (0.623489801858733, -0.78183148246803) circle (2.2630899352994pt);
\fill[black] (2, 0) circle (2.2630899352994pt);
\fill[black] (0.623489801858734, 0.78183148246803) circle (2.2630899352994pt);
\fill[black] (0.900968867902419, 0.433883739117558) circle (2.2630899352994pt);
\fill[black] (-0.222520933956315, -0.974927912181824) circle (2.2630899352994pt);
\fill[black] (-0.900968867902419, -0.433883739117558) circle (2.2630899352994pt);
\fill[black] (-0.222520933956314, 0.974927912181824) circle (2.2630899352994pt);
\fill[black] (0.900968867902419, -0.433883739117558) circle (2.2630899352994pt);
\fill[black] (0.222520933956314, 0.974927912181824) circle (2.2630899352994pt);
\fill[black] (0.222520933956313, -0.974927912181824) circle (2.2630899352994pt);
\fill[black] (-0.900968867902419, 0.433883739117558) circle (2.2630899352994pt);
\fill[black] (-1.24697960371747, 1.56366296493606) circle (2.2630899352994pt);
\fill[black] (-1, 1.22464679914735e-16) circle (2.2630899352994pt);
\end{tikzpicture}

%% file: 15.tex
\begin{tikzpicture}[scale=0.5]
\draw (0,1) -- (0.4067366430758,0.913545457642601);
\draw (0,1) -- (-0.4067366430758,0.913545457642601);
\draw (0.4067366430758,0.913545457642601) -- (0.743144825477394,0.669130606358858);
\draw (-0.4067366430758,0.913545457642601) -- (-0.743144825477394,0.669130606358858);
\draw (-0.587785252292473,-0.809016994374948) -- (-0.207911690817759,-0.978147600733806);
\draw (-0.587785252292473,-0.809016994374948) -- (-0.866025403784438,-0.5);
\draw (-0.207911690817759,-0.978147600733806) -- (0.207911690817759,-0.978147600733806);
\draw (-0.866025403784438,-0.5) -- (-0.994521895368273,-0.104528463267654);
\draw (-0.743144825477394,0.669130606358858) -- (-0.951056516295154,0.309016994374947);
\draw (-0.951056516295154,0.309016994374947) -- (-0.994521895368273,-0.104528463267654);
\draw (-0.951056516295154,0.309016994374947) -- (-1.90211303259031,0.618033988749894);
\draw (0.743144825477394,0.669130606358858) -- (0.951056516295154,0.309016994374947);
\draw (0.207911690817759,-0.978147600733806) -- (0.587785252292473,-0.809016994374947);
\draw (0.951056516295154,0.309016994374947) -- (0.994521895368273,-0.104528463267653);
\draw (0.587785252292473,-0.809016994374947) -- (0.866025403784439,-0.5);
\draw (0.994521895368273,-0.104528463267653) -- (0.866025403784439,-0.5);
\fill[black] (-0.866025403784438, -0.5) circle (2.12809111797729pt);
\fill[black] (0.587785252292473, -0.809016994374947) circle (2.12809111797729pt);
\fill[black] (0.207911690817759, -0.978147600733806) circle (2.12809111797729pt);
\fill[black] (-0.951056516295154, 0.309016994374947) circle (2.12809111797729pt);
\fill[black] (0.994521895368273, -0.104528463267653) circle (2.12809111797729pt);
\fill[black] (-0.587785252292473, -0.809016994374948) circle (2.12809111797729pt);
\fill[black] (0.866025403784439, -0.5) circle (2.12809111797729pt);
\fill[black] (-0.207911690817759, -0.978147600733806) circle (2.12809111797729pt);
\fill[black] (-0.994521895368273, -0.104528463267654) circle (2.12809111797729pt);
\fill[black] (0.4067366430758, 0.913545457642601) circle (2.12809111797729pt);
\fill[black] (-0.4067366430758, 0.913545457642601) circle (2.12809111797729pt);
\fill[black] (0.743144825477394, 0.669130606358858) circle (2.12809111797729pt);
\fill[black] (-1.90211303259031, 0.618033988749894) circle (2.12809111797729pt);
\fill[black] (0, 1) circle (2.12809111797729pt);
\fill[black] (0.951056516295154, 0.309016994374947) circle (2.12809111797729pt);
\fill[black] (-0.743144825477394, 0.669130606358858) circle (2.12809111797729pt);
\end{tikzpicture}

%% file: 16.tex
\begin{tikzpicture}[scale=0.5]
\draw (0,1) -- (0.38268343236509,0.923879532511287);
\draw (0,1) -- (-0.38268343236509,0.923879532511287);
\draw (0.38268343236509,0.923879532511287) -- (0.707106781186547,0.707106781186548);
\draw (-0.38268343236509,0.923879532511287) -- (-0.707106781186548,0.707106781186547);
\draw (-0.38268343236509,-0.923879532511287) -- (1.22464679914735e-16,-1);
\draw (-0.38268343236509,-0.923879532511287) -- (-0.707106781186547,-0.707106781186548);
\draw (1.22464679914735e-16,-1) -- (0.38268343236509,-0.923879532511287);
\draw (-0.707106781186547,-0.707106781186548) -- (-0.923879532511287,-0.38268343236509);
\draw (-0.707106781186548,0.707106781186547) -- (-0.923879532511287,0.38268343236509);
\draw (-1,-1.83697019872103e-16) -- (-0.923879532511287,-0.38268343236509);
\draw (-1,-1.83697019872103e-16) -- (-0.923879532511287,0.38268343236509);
\draw (0.707106781186547,0.707106781186548) -- (0.923879532511287,0.38268343236509);
\draw (0.38268343236509,-0.923879532511287) -- (0.707106781186548,-0.707106781186547);
\draw (0.923879532511287,0.38268343236509) -- (1,6.12323399573677e-17);
\draw (0.707106781186548,-0.707106781186547) -- (0.923879532511287,-0.38268343236509);
\draw (1,6.12323399573677e-17) -- (0.923879532511287,-0.38268343236509);
\fill[black] (-0.707106781186548, 0.707106781186547) circle (2pt);
\fill[black] (0.923879532511287, 0.38268343236509) circle (2pt);
\fill[black] (-0.923879532511287, 0.38268343236509) circle (2pt);
\fill[black] (0, 1) circle (2pt);
\fill[black] (0.707106781186547, 0.707106781186548) circle (2pt);
\fill[black] (-0.38268343236509, 0.923879532511287) circle (2pt);
\fill[black] (0.38268343236509, 0.923879532511287) circle (2pt);
\fill[black] (1.22464679914735e-16, -1) circle (2pt);
\fill[black] (-0.923879532511287, -0.38268343236509) circle (2pt);
\fill[black] (0.923879532511287, -0.38268343236509) circle (2pt);
\fill[black] (-0.38268343236509, -0.923879532511287) circle (2pt);
\fill[black] (-1, -1.83697019872103e-16) circle (2pt);
\fill[black] (1, 6.12323399573677e-17) circle (2pt);
\fill[black] (0.38268343236509, -0.923879532511287) circle (2pt);
\fill[black] (-0.707106781186547, -0.707106781186548) circle (2pt);
\fill[black] (0.707106781186548, -0.707106781186547) circle (2pt);
\end{tikzpicture}